\documentclass[a4paper,reqno,11pt]{amsart}
\usepackage[utf8]{inputenc}
\usepackage[T1]{fontenc}
\usepackage{amsfonts}
\usepackage{graphicx}
\usepackage{amsmath}
\usepackage{caption}
\usepackage[normalem]{ulem} 
\usepackage{subcaption}
\usepackage{float}
\usepackage{amsmath}
\numberwithin{equation}{section}
\usepackage{xcolor}
\usepackage{amsmath}
\usepackage{tikz}
\usepackage{comment}
 \usepackage[foot]{amsaddr}
\usepackage{epsfig}
\usepackage{color,hyperref}
\usetikzlibrary{fit,
  matrix,
  positioning,
  decorations,
  decorations.pathreplacing
}

\usepackage{latexsym}
\definecolor{MyOra}{RGB}{255, 165, 0}
\definecolor{MyOra1}{RGB}{255,119,0}
\definecolor{MyBlue}{rgb}{0.25,0.5,0.75}

\newcommand{\myrmk}[2]{\begin{quote}\fbox{\begin{minipage}{\linewidth}\small \noindent{\color{red}\bf #2:} #1{\color{red}}\end{minipage}}\end{quote}}
\newcommand{\rmkH}[1]{\myrmk{#1}{HD}}
\newcommand{\rmkP}[1]{\myrmk{#1}{PR}}

\newcounter{rhp}
\newenvironment{rhp}[1][]{\refstepcounter{rhp}\par\medskip
   \noindent \textbf{Riemann-Hilbert problem~\therhp. #1} \rmfamily}{\medskip}

\makeatletter

\newcommand{\tr}{\operatorname{tr}}
\usepackage{xparse}
\DeclareUnicodeCharacter{00A0}{~}
\usepackage{hyperref}
\usepackage[margin=1in]{geometry}
\usepackage{bbold}
\usepackage[firstinits=true,uniquename=false,uniquelist=false,style=alphabetic,sortcites=false,natbib=true,doi=true,backend=biber]{biblatex}
\title{On a class of elliptic orthogonal polynomials and their integrability}

\author{Harini Desiraju$^1$}
 \email{\href{mailto:harini.desiraju@sydney.edu.au}{harini.desiraju@sydney.edu.au}}
\author{Tomas Lasic Latimer$^1$}
\email{\href{mailto:mailto:t.lasiclatimer@maths.usyd.edu.au}{t.lasiclatimer@maths.usyd.edu.au}}
\author{Pieter Roffelsen$^1$}
 \email{\href{mailto:pieter.roffelsen@sydney.edu.au}{pieter.roffelsen@sydney.edu.au}}
\address{$^1$School of Mathematics and Statistics, University of Sydney, Camperdown, NSW 2006}

\date{}

\newtheorem{lemma}{Lemma}[section]

\newtheorem{theorem}{Theorem}[section]
\newtheorem{corollary}{Corollary}[section]
\newtheorem{proposition}{Proposition}[section]

\newtheorem{remark}{Remark}[section]

\begin{document}

\maketitle
\begin{abstract}
Building upon the recent works of
Bertola; Fasondini, Olver and Xu, we define a class of orthogonal polynomials on elliptic curves and establish a corresponding Riemann-Hilbert framework.
We then focus on the special case, defined by a constant weight function, and use the Riemann-Hilbert problem to derive recurrence relations and differential equations for the orthogonal polynomials.
We further show that the sub-class of even polynomials is associated to the elliptic form of Painlev\'e VI, with the tau function given by the Hankel determinant of even moments, up to a scaling factor.
The first iteration of these even polynomials relates to the special case of Painlev\'e VI studied by Hitchin in relation to self-dual Einstein metrics.

 \end{abstract}

\tableofcontents
\pagebreak

\section{Introduction}
Orthogonal polynomials constitute a fundamental class of special functions with important applications to a wide array of topics, from combinatorics to signal processing. 
Particularly, they provide useful tools to understand universality of random matrix ensembles \cite{deift1999orthogonal} and large $N$ limits of matrix models \cite{di19952d},
and describe special solutions to integrable systems such as Painlev\'e equations \cite{van2017orthogonal}.


Traditionally, orthogonal polynomials define a basis of real polynomials, orthogonal with respect to an inner product defined by integrating against a weight function on a subset of the real line. For example, Hermite polynomials are orthogonal with respect to the weight $e^{-x^2}$ on the real line. 
Moreover, it is known that generalisations of classical Chebyshev and Jacobi polynomials can be described by elliptic functions, see for instance
\cite{akhiezer1990elements, carlitz1960some,ismail2001some, rees1945elliptic, vinet2009elliptic}.

In the past couple of years there were notable breakthroughs in defining orthogonal polynomials directly on elliptic curves,
with modern techniques facilitating a systematic analysis of their properties \cite{bertolaPade2021,bertola2022nonlinear,bertola2022critical,fasondini2022orthogonal,fasondini2023orthogonal}. Inspired by these works, we consider families of meromorphic functions built out of the the Weierstrass $\wp$-function and its derivative, which are orthogonal with respect to a given weight function, and call them elliptic orthogonal polynomials (EOPs) due to their proximity in construction of those introduced by Heine \cite{heine1878handbuch} and Rees \cite{rees1945elliptic}. 

We establish a general framework to analyse such polynomials using their moments and the Riemann-Hilbert method.
Furthermore, when the weight is constant, we show that the even EOPs, indexed by $k$, are related to the elliptic form of Painlev\'e VI. For $k=1$, the parameters of the elliptic Painlev\'e VI equation are $\left( \frac{1}{8}, \frac{1}{8}, -\frac{1}{8}, \frac{3}{8}\right)$, for which the general solution is known to be described by elliptic functions \cite{hitchin1995twistor, manin1998universal}.


The notion of EOPs we use in this paper is as follows: let $\tau$ be an element of the upper half-plane $\mathbb{H}$ and $\pi(z)$ be an elliptic function with periods $1$ and $\tau$. We call $\pi(z)$ an elliptic polynomial if all of its poles are located on the lattice $\mathbb{Z}+\mathbb{Z}\cdot \tau$.
Its degree is $n$ if the pole at $z=0$ is of order $n$, and we call it monic if, for $n\geq 0$,
\begin{equation*}
    \pi(z)=z^{-n}(1+\mathcal{O}(z)), \qquad \textrm{as} \quad z\to 0.
\end{equation*}
 Note, in particular, that there exists no elliptic polynomial of degree one.
We consider sequences of elliptic polynomials $(\pi_n)_{n\geq 0,n\neq 1}$, with $\pi_n$ monic of degree $n$ for $n\in\mathbb{N}_{\neq 1}$, which satisfy an orthogonality condition of the form
\begin{equation}\label{orthogonality}
    \int_{\frac{\tau}{2}}^{\frac{\tau}{2}+1}\pi_m(z)\pi_n(z){\sf w}(z)dz=\delta_{mn}h_n,
\end{equation}
for some $h_n\in\mathbb{C}$, where $\delta_{mn}$ is the Kronecker delta function and $w(z)$ is an $L^1$ function on the interval $\gamma:=[\tfrac{\tau}{2},\tfrac{\tau}{2}+1]$, called the weight function, for all $m,n\in\mathbb{N}_{\neq 1}$. In this case, we call $(\pi_n)_{n\geq 0,n\neq 1}$ a sequence of elliptic orthogonal polynomials (EOPs).
The choice of support $\gamma$ is motivated by the fact that it is invariant under negation and complex conjugation on the torus
\begin{equation*}
    \mathbb{T}:=\mathbb{C}/(\mathbb{Z}+\mathbb{Z}\cdot \tau),
\end{equation*}
and does not contain $[0]\in\mathbb{T}$.

We note the following analogy between elliptic polynomials and traditional complex polynomials. A complex polynomial can be characterised as a meromorphic function on $\mathbb{CP}^1$, with at most one pole, at $\infty$. Analogously, an elliptic polynomial, as defined above, can be characterised as a meromorphic function on the torus $\mathbb{T}$, with only at most one pole, at $[0]\in \mathbb{T}$.

 \begin{figure}[h]
\includegraphics[trim={0cm 0cm 0cm 0cm},clip, width=7cm]{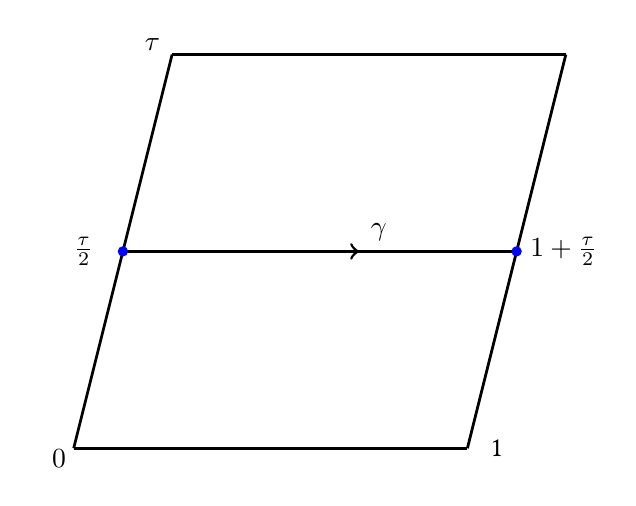}
    \caption{The orthogonality interval in the fundamental domain $\mathbb{T}$.}
    \label{fig:1}
\end{figure}


A useful basis for elliptic polynomials can be constructed in terms of the Weierstrass $\wp$-function  and its $z$-derivative,
\begin{align}\label{def:basis}
\mathcal{B}=\{\mathcal{E}_n\}_{{n\geq 0,n\neq 1}}, &&\mathcal{E}_{2k}=\wp(z)^k, && \mathcal{E}_{2k+3}=-\tfrac{1}{2}\wp'(z)\wp(z)^k , && k\geq 0,
\end{align}
chosen such that $\mathcal{E}_n$ is monic of degree $n$, for $n\in\mathbb{N}_{\neq 1}$. 
Recall the isomorphism
\begin{align}\label{transf:xy}
x =\wp(z),&&y = \wp'(z),
\end{align}
from $\mathbb{T}$ to the cubic curve
   \begin{equation}\label{OG elliptic curve}
        y^2 = 4x^3 - g_2x - g_3,
    \end{equation} 
where $g_{2,3}$ are the elliptic invariants. Under this map,  the EOP $\pi_n$ becomes a bivariate polynomial in $\{x,y\}$ of weighted degree $n$, where
\begin{equation*}
    \operatorname{deg}(x)=2,\quad \operatorname{deg}(y)=3.
\end{equation*}
Under this identification, our definition of elliptic polynomials then coincides with that used in \cite{fasondini2022orthogonal},
and the transformation \eqref{transf:xy} results in an equivalent definition of the bases of elliptic polynomials to that given in \eqref{def:basis}.

When ${\sf w}(z)$ is even in $z$ around the midpoint of the contour $\gamma$, that is,
\begin{equation}\label{eq:even_weight}
    {\sf w}(\tfrac{1}{2}(1+\tau)+z)={\sf w}(\tfrac{1}{2}(1+\tau)-z)\quad z\in [0,1],
\end{equation}
the monic elliptic-polynomials 
naturally split into even and odd polynomials:
\begin{align}
    \pi_{2k}(z,\tau) &= \sum_{i=0}^k \widehat{a}_{i,2k}(\tau) \wp(z)^{k-i}, \qquad \widehat{a}_{0,2k}=1,\\
    \pi_{2k+3}(z,\tau) &= -\sum_{i=0}^k \widehat{a}_{i,2k+3}(\tau) \wp(z)' \wp(z)^{k-i}/2, \qquad \widehat{a}_{0,2k+3}=1,
\end{align}
where,
\begin{align}\label{oddsandevens}
    \pi_{2k}(-z) = \pi_{2k}(z), && \pi_{2k+3}(-z) = - \pi_{2k+3}(z).
\end{align}
Henceforth, we refer to $\pi_{2k}$ and $\pi_{2k+3}$ as the respective even and odd EOPs. Note that $\pi_{2k}(z)$ is purely a function of $x=\wp(z)$ according to \eqref{def:basis}, and such polynomials are similar in nature to the Akheizer polynomials
 in \cite{CHENIT2008}, and generalised Jacobi polynomials studied in the literature by \cite{carlitz1960some,heine1878handbuch,ismail2001some,rees1945elliptic} among others.

In this paper, we show how the EOPs, for any choice of weight function, can be written in terms of determinants of moments and characterised as the $(1,1)$ entry of the unique solution of a corresponding $2\times 2$ Riemann-Hilbert problem (RHP). 

The connection between orthogonal polynomials and the RHP was established in the 1990's by Fokas, Its, and Kitaev \cite{fokas1991discrete}. Since then it has been instrumental in a number of settings to prove a variety of results for different classes of orthogonal polynomials \cites{Deift1999strong,qRHP,etnavol25pp369-392}, most often finding application in determining the asymptotics of different classes of orthogonal polynomials, and universality of random matrix ensembles.
RHPs were also used in \cite{bertola2022nonlinear} to determine large degree asymptotic behaviours within an analytically distinct class of orthogonal polynomials on elliptic curves, allowing for an additional pole in the basis. 

A family of orthogonal polynomials corresponding to a given weight function are completely determined by the moments of the weight function \cite{szeg1939orthogonal}. Such a representation plays an important role in several aspects of orthogonal polynomials, as well as their applications in areas such as combinatorics and probability theory, see the survey articles \cite{corteel2016moments, konig2005orthogonal} for example. Specifically, the Hankel determinants of moments were exploited to study the partition functions of ensembles of random matrices and to study the corresponding tau-functions or, equivalently, the solutions of integrable equations. See Forrester and Witte \cite{forrester2006random} for this approach in the case of Painlev\'e VI and Bertola \cite{bertola2008moment} for other classes of isomonodromic systems.
In the present case, the relevant moments and the determinant are defined as 
\begin{align*}
    \mu_{i,j} := \int_{\gamma} \mathcal{E}_i(z)\mathcal{E}_j(z) {\sf w}(z)dz, && D_n:=\det(\mu_{i,j})_{i,j=0}^{n-1}.
\end{align*}
As is the case for classical orthogonal polynomials \cite{fokas1991discrete}, the aforementioned RHP for EOPs is uniquely solvable if and only if the determinant $D_n \neq 0$. Furthermore, the moment matrix is a block matrix consisting of a checkerboard pattern of odd and even moments. Consequently its determinant factorises into two Hankel determinants, consisting of even and odd moments respectively, both of which define Painlev\'e tau-functions.  
 A similar result connecting Hankel determinants constructed out of certain elliptic functions and the tau-function of Painlev\'e VI was obtained by \cite{basor2015asymptotics} from the study of generalized Jacobi polynomials.

\subsection{Outlook}
 We list a few possible future directions following the results presented here:
\begin{enumerate}
    \item with the construction in Section \ref{Sec:RHP_HD}, one can systematically increase the complexity of the weight function. For example, 
    a weight function that is algebraic in $\wp(z)$, which may lead to solutions of other integrable equations,
    \item the large degree asymptotics of EOPs can be studied using the RHPs described here.
    \item  Hankel determinants are related to partition functions of random matrix ensembles. It would be interesting to see the elliptic extensions of such relations.
    \item Another natural question would be the extension of the Riemann-Hilbert setup developed here to study orthogonal polynomials on higher genus surfaces \cite{fasondini2022orthogonal}.
\end{enumerate}

\subsection{Outline}
In Section \ref{Sec:RHP_HD} we begin by describing the general and even polynomials in terms of the determinants of their moments, and the construct the solutions $Y_n$ and $Y_{2k}$ of the Riemann-Hilbert problems  associated to each class of polynomials respectively. Theorem \ref{thm:existence_uniqueness} proves the existence and uniqueness of the solution of the general RHP. For the remainder of the paper we set the weight ${\sf w}(z)\equiv 1$. 

In Section \ref{Sec:Mix_Lin}, we obtain  differential and discrete linear systems satisfied by the solution of the general RHP $Y_n$ in Propositions \ref{prop:difference_ls_Yn} and \ref{prop:differential_ls_Yn} respectively. In Theorems  \ref{thm:rec_pol_gen} and \ref{main theorem sect 3} we derive recurrence relations for the polynomials $\pi_n$ and the  coefficients.

In Section \ref{Sec:Even_Lin}, Theorem \ref{thm:Lax_ell_PVI} shows that the linear problems satisfied by the solution of the even RHP $Y_{2k}$ describe the Lax pair of the elliptic form of Painlev\'e VI. For $k=1$, the monodromy exponents assume a special form which is related to the Hitchin case of Painlev\'e VI as shown in Proposition \ref{prop:PVIell_normal}.

In Section \ref{Sec:HD_PVI}, we derive explicit solutions of Painlev\'e VI built out of Hankel determinants of even moments, see Theorem \ref{thm:pvisolution}. Furthermore, we show that these Hankel determinants are the corresponding Painlev\'e VI tau-functions, in Theorem \ref{thm:tau}.

\subsection{Acknowledgements}
We thank Nalini Joshi and Milena Radnovi\'c for valuable discussions and suggestions.
H.D. acknowledges the support of Australian Research Council Discovery Project \#DP200100210. P.R. acknowledges the support of Australian Research
Council Discovery Project \#DP210100129. TL's research was supported the Australian Government Research Training Program and by the University of Sydney Postgraduate Research Supplementary Scholarship in Integrable Systems.

A part of the work was done during the authors' residence at the Isaac Newton Institute during the Fall 2022 semester and they thank the Møller institute and organisers of the program ”Applicable resurgent asymptotics: towards a universal theory” for their hospitality. HD thanks the Simons foundation for supporting her INI visit.


\section{Riemann-Hilbert problems and moments}\label{Sec:RHP_HD}


In this section, we begin by detailing the representation of the general and even EOPs in terms of determinants of moments.
We then show show that the consecutive EOPs, and their suitable Cauchy transforms, form a unique solution to a corresponding Riemann-Hilbert problem (RHP). We relate the existence and uniqueness of the polynomials to the unique solvability of the RHP and find explicit expressions for the determinants of the solutions. This will in turn allow us to derive differential and difference linear systems satisfied by the polynomials for the case ${\sf w}(z)\equiv 1$.


\subsection{Moments in the general case}
We start by assuming $\tau\in i \mathbb{R}$, so that the Weierstrass $\wp$-function is real on $\gamma$, {\it i.e} the basis $\mathcal{B}$ of elliptic polynomials consists of real functions on $\gamma$. If, in addition, the weight function ${\sf w}(z)$ is strictly positive, then
\begin{equation*}
    \langle f,g\rangle=\int_\gamma{f(z)g(z){\sf w}(z)dz}
\end{equation*}
defines an inner product on the space of real elliptic polynomials and the Gram–Schmidt process shows that the corresponding EOPs exist and are unique, with $\pi_n$ given by
\begin{equation}
 \label{OG delta}
\pi_{n}(z) = D_n^{-1}
\begin{vmatrix}
\mu_{0,0} & \mu_{0,2} & \mu_{0,3} &\hdots&\mu_{0,n} \\
\mu_{2,0} & \mu_{2,2} & \mu_{2,3}&\hdots&\mu_{2,n} \\
\mu_{3,0} & \mu_{3,2} & \mu_{3,3}&\hdots&\mu_{3,n} \\
\vdots & \vdots &\vdots & \ddots &\vdots \\
\mu_{n-1,0} & \mu_{n-1,1} & \mu_{n-1,2} &\hdots&\mu_{n-1,n} \\
\mathcal{E}_0(z) & \mathcal{E}_2(z)  &\mathcal{E}_3(z)  &...& \mathcal{E}_{n}(z) 
\end{vmatrix},
\end{equation}
where
\begin{equation}\label{def:moment}
    \mu_{i,j}=\int_\gamma \mathcal{E}_i(z)\mathcal{E}_j(z) {\sf w}(z) dz\qquad (i,j\in\mathbb{N}_{\neq 1}),
\end{equation}
and 
\begin{equation}\label{eq:Dn}
D_n=|S_n|,\quad S_n:= \begin{vmatrix}
\mu_{0,0} & \mu_{0,2} & \mu_{0,3} &\hdots&\mu_{0,n-1} \\
\mu_{2,0} & \mu_{2,2} & \mu_{2,3}&\hdots&\mu_{2,n-1} \\
\mu_{3,0} & \mu_{3,2} & \mu_{3,3}&\hdots&\mu_{3,n-1} \\
\vdots & \vdots &\vdots & \ddots &\vdots \\
\mu_{n-1,0} & \mu_{n-1,1} & \mu_{n-1,2} &\hdots&\mu_{n-1,n-1} 
\end{vmatrix}.
\end{equation}
See Appendix \ref{appendix:mom} for an illustration.
The following identity then follows from the above equation,
\begin{equation}\label{eq:normDn}
h_{n}:=\int_{\frac{1}{2}\tau+0}^{\frac{1}{2}\tau+1}\pi_{n}(z)^2{\sf w}(z)dz=\int_{\frac{1}{2}\tau+0}^{\frac{1}{2}\tau+1}\mathcal{E}_n\pi_{n}(z){\sf w}(z)dz=\frac{D_{n+1}}{D_{n}}.
\end{equation}
For $\tau \notin i \mathbb{R}$, or ${\sf w}(z)$ is not strictly positive, the polynomial $\pi_n(z)$ exists and is unique if and only if the determinant $D_n$
is nonzero. Similarly, we will find that the associated RHP has a unique solution if and only if $D_n\neq 0$.

\begin{figure}[h]
    \centering
    \includegraphics[width=0.7\textwidth]{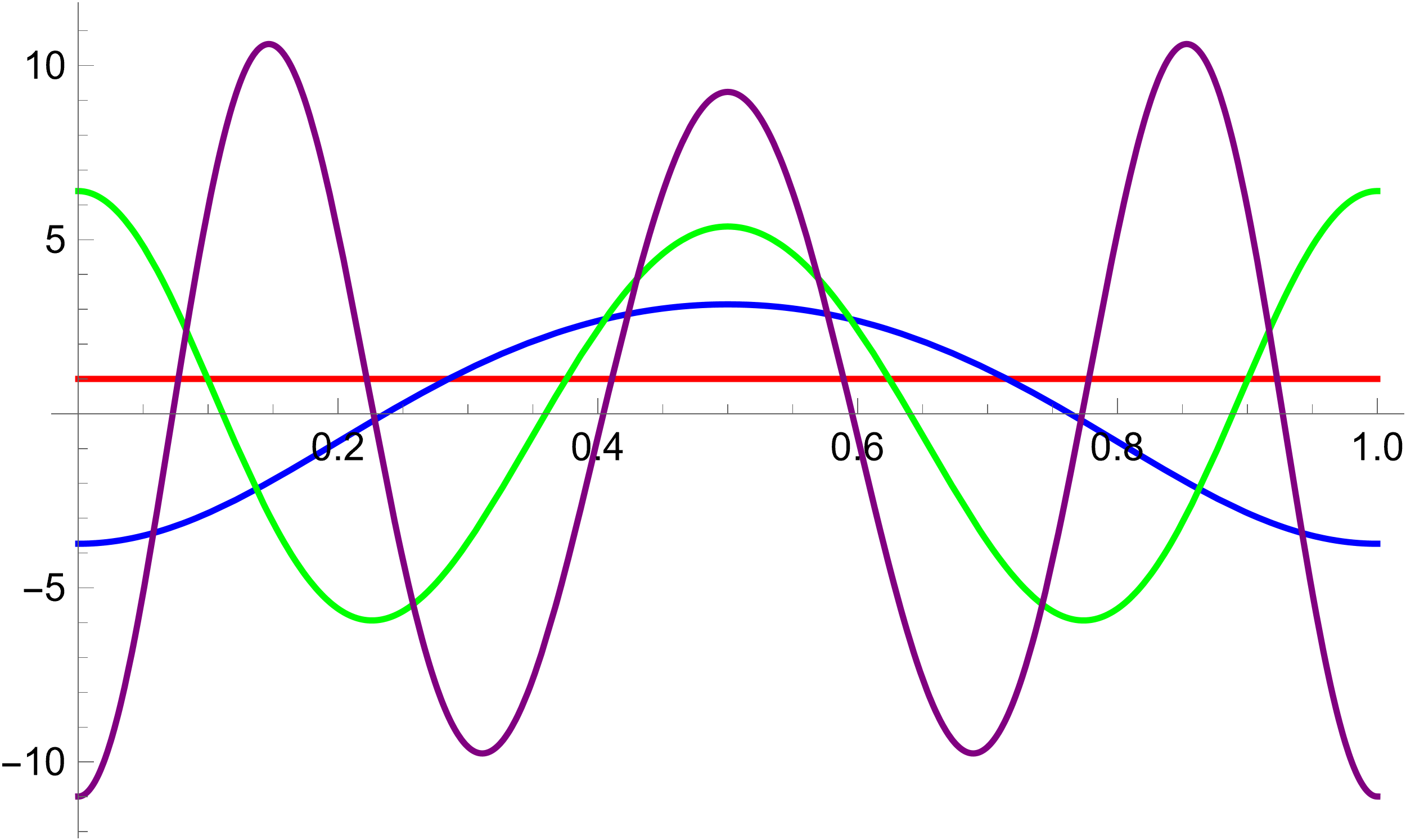}
    \caption{Graph displaying function values of the even EOPs $\pi_0$, $\pi_2$, $\pi_4$ and $\pi_6$, on the interval $\tfrac{1}{2}\tau+[0,1]$, in red, blue, green and purple respectively, with $\tau=i$.}
    \label{fig:even}
\end{figure}

\begin{figure}[h]
    \centering
    \includegraphics[width=0.7\textwidth]{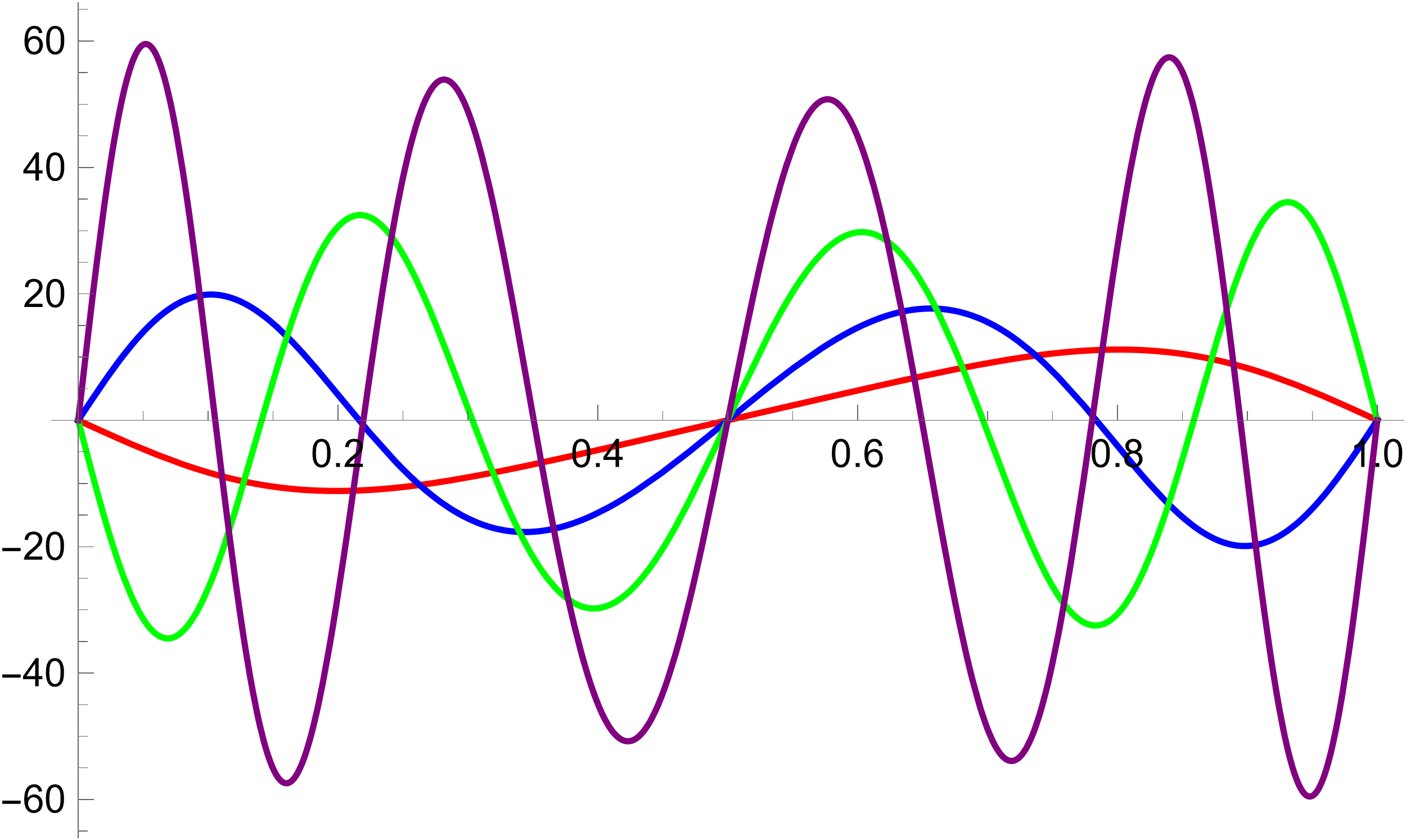}
    \caption{Graph displaying function values of the odd EOPs $\pi_3$, $\pi_5$, $\pi_7$ and $\pi_9$, on the interval $\tfrac{1}{2}\tau+[0,1]$, in red, blue, green and purple respectively, with $\tau=i$.}
    \label{fig:odd}
\end{figure}

\begin{figure}[H]
    \centering
    \includegraphics[width=0.7\textwidth]{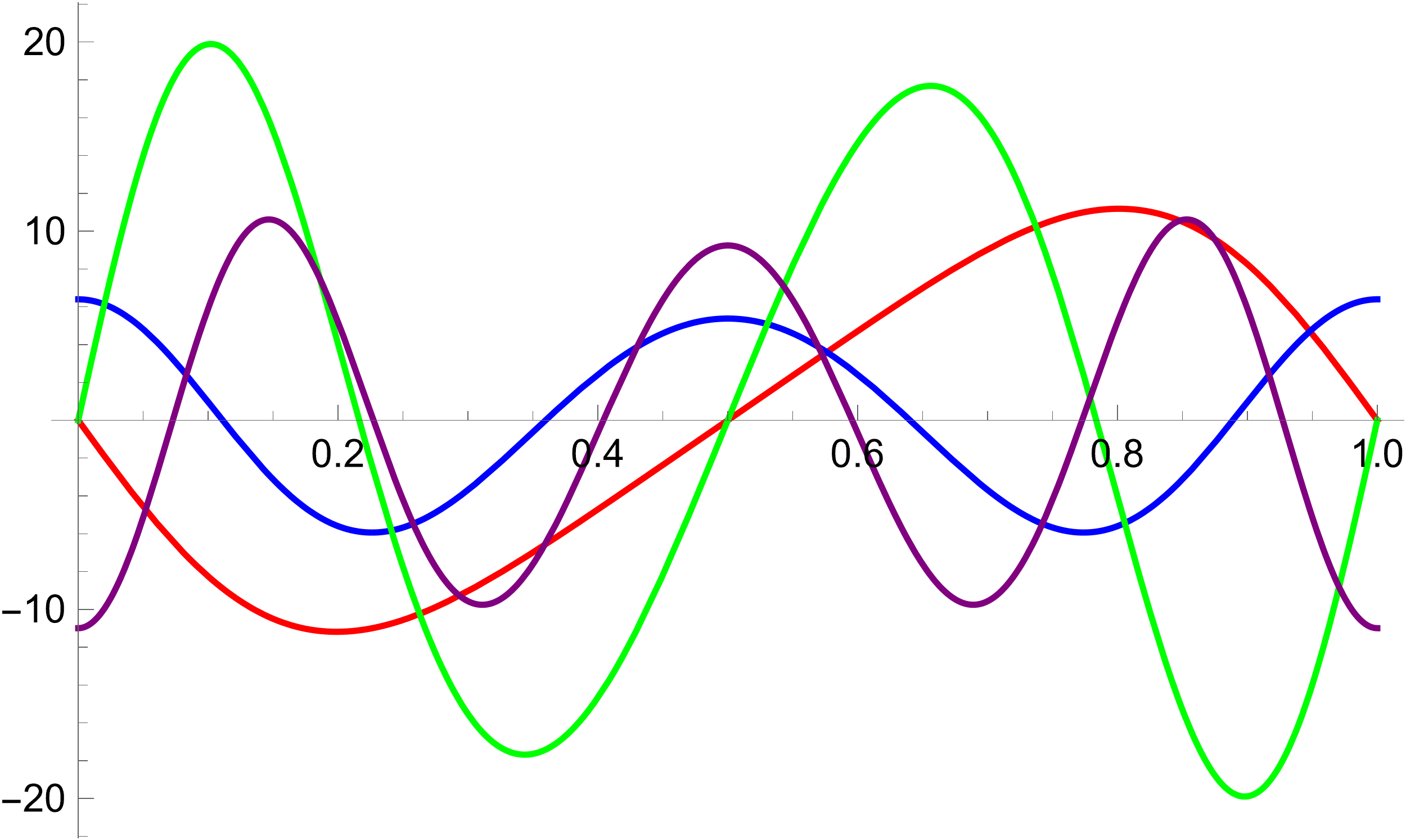}
    \caption{Graph displaying function values of the EOPs $\pi_3$, $\pi_4$, $\pi_5$ and $\pi_6$, on the interval $\tfrac{1}{2}\tau+[0,1]$, in red, blue, green and purple respectively, with $\tau=i$.}
    \label{fig:mixed}
\end{figure}

In Figures \ref{fig:even}, \ref{fig:odd} and \ref{fig:mixed}, the first couple of EOP's are plotted on $\gamma$, with $\tau=i$ and ${\sf w}(z)\equiv 1$.


\subsection{Riemann-Hilbert problem for the general case}
The first step in defining the RHP is to analyse the Cauchy transform of the polynomials 
\begin{align}\label{def:cauchypin}
     \mathcal{C}(\pi_n)(z):= \int_{\gamma} \pi_n(w) C(w,z) {\sf w}(w)\frac{dw}{2\pi i}.
\end{align}
Generally, the Cauchy kernel is a meromorphic function in $z$ and a one form in $w$, with residues $\pm 1$ at its poles, {\it i.e} it is an Abelian differential of the third kind. The precise form of such a kernel on a torus is not unique and several examples can be found in the literature \cite{bikbaev1988asymptotics, rodin1987riemann,del2020isomonodromic,bertola2022nonlinear}. For the present case, we consider the scalar kernel 
\begin{gather}\label{def:cauchyker}
    C(w,z) = \zeta(w-z)-\zeta(w),
\end{gather}
where $\zeta(.)$ denotes the Weierstrass $\zeta$-function (see \ref{Appendix:EllFunc}). This Cauchy kernel has the following properties.
\begin{enumerate}
\item The periodicity of $\zeta$-function, see equation \eqref{periodicity:zeta}, implies that 
\begin{align}\label{periodicity:Cauchy}
    C(w, z+1) = C(w,z) - \eta_{1}(\tau),&&  C(w, z+\tau) = C(w, z) - \eta_{2}(\tau).
\end{align}
\item The kernel $C(w,z)$ has poles at $w=z, w=0$ with residues $\pm 1$ respectively, 
\item In the limit $w\to z$,
\begin{gather}
   C(w,z) = \frac{1}{w-z} +\mathcal{O}(w-z), 
\end{gather}
\item and in the limit $z\to 0$, 
\begin{align}\label{eq:cauchytaylor}
    C(w,z) = \sum_{k=1}^{\infty} (-1)^{k} \frac{z^k}{k!} \zeta^{(k)}(w),
\end{align}
where $\zeta^{(k)}(\cdot)$ denotes the $k$th derivative of $\zeta$.
\end{enumerate}
The following statements then  hold true for \eqref{def:cauchypin}.

\begin{lemma}\label{lemma:periodicity}
    The Cauchy transform of $\pi_n$ is a piecewise analytic, doubly periodic function for $n\geq 2$.
\end{lemma}
\begin{proof} It is immediate from the definition, that $\mathcal{C}(\pi_n)(z)$ is analytic away from $\gamma$.
From the relations \eqref{periodicity:Cauchy}, we see that
\begin{align}
    \mathcal{C}(\pi_n)(z+1) = \mathcal{C}(\pi_n)(z) - \int_{\gamma} \eta_1(\tau) \pi_{n}(w) {\sf w}(w)\frac{dw}{2\pi i},
\end{align}
and due to orthogonality \eqref{orthogonality}, the second term vanishes for $n\neq 0$. Similar computation holds for $\mathcal{C}(\pi_n)(z+\tau)$.
\end{proof}

\begin{lemma}\label{Series lemma}
    In the asymptotic limit $z\to 0$, the polynomials have an asymptotic expansion of the form
    \begin{align}\label{asymp:pin}
        \pi_n(z)= \sum_{j=0}^{\infty} \frac{c_{j,n}(\tau)}{z^{n-j}}, \qquad c_{0,n} =1, 
    \end{align}
     and, similarly, their Cauchy transform have an asymptotic expansion of the form,\footnote{all the coefficients are $\tau$-dependent unless stated otherwise. We often omit the $\tau$- dependence of the coefficients for ease of notation.}
    \begin{align}\label{asymp:Cpinmixed}
        \mathcal{C}(\pi_n)(z)=  \frac{h_n(\tau)}{2\pi i} \sum_{j=0}^{\infty} \widetilde{c}_{j,n}(\tau) z^{n+j-1}, && \widetilde{c}_{0,n} =1, &&\widetilde{c}_{1,n} =-c_{1, n+1}.
    \end{align}
\end{lemma}
\begin{proof}
Using the Taylor expansion around $z=0$ of the Cauchy kernel in equation \eqref{eq:cauchytaylor}, which holds uniformly in $w\in\gamma$, we obtain the following asymptotic expansion as $z\rightarrow 0$,
\begin{align}
\mathcal{C}(\pi_n)(z) &=\frac{1}{2\pi i} \int_{\gamma}    \sum_{k=1}^{\infty} (-1)^{k} \frac{z^k}{k!} \zeta^{(k)}(w)\pi_n(w) {\sf w}(w) \frac{dw}{2\pi i} \\
&=\frac{1}{2\pi i}    \sum_{k=1}^{\infty} (-1)^{k} \frac{z^k}{k!}\int_{\gamma} \zeta^{(k)}(w)\pi_n(w) {\sf w}(w) \frac{dw}{2\pi i} \\
&=\frac{1}{2\pi i}    \sum_{k=n-1}^{\infty} (-1)^{k} \frac{z^k}{k!}\int_{\gamma} \zeta^{(k)}(w)\pi_n(w) {\sf w}(w) \frac{dw}{2\pi i} \\
& = \frac{h_n(\tau)}{2\pi i} \sum_{j=0}^{\infty} \widetilde{c}_{j,n}(\tau) z^{n+j-1},
\end{align}
where, in the third equality we used orthogonality and the fact that $\zeta^{(k)}(w)$ is an elliptic polynomial of degree $k$, for $k\geq 1$, with coefficients given by
\begin{equation*}
   \frac{h_n(\tau)}{2\pi i}\widetilde{c}_{j,n}(\tau)=\frac{(-1)^{n-1+j}}{(n-1+j)!} \int_{\gamma} \zeta^{(n-1+j)}(w)\pi_n(w) {\sf w}(w) \frac{dw}{2\pi i},
\end{equation*}
for $j\geq 0$.
To compute the leading order coefficients, note that
\begin{equation*}
    \frac{(-1)^{n-1+j}}{(n-1+j)!} \zeta^{(n-1+j)}(w)=w^{-(n+j)}(1+\mathcal{O}(w^2)),
\end{equation*}
 as $w\rightarrow 0$, so that
 \begin{equation*}
    \frac{(-1)^{n-1}}{(n-1)!} \int_\gamma \zeta^{(n-1)}(w)\pi_n(w){\sf w}(w)dw=\int_\gamma \pi_n(w)^2dw=h_n,
 \end{equation*}
 and
  \begin{equation*}
    \frac{(-1)^{n}}{n!} \int_\gamma \zeta^{(n)}(w)\pi_n(w){\sf w}(w)dw=\int_\gamma (\pi_{n+1}(w)-c_{1,n+1}\pi_n(z))\pi_n(w){\sf w}(w)dw=-c_{1,n+1}h_n.
 \end{equation*}
Therefore, $\widetilde{c}_{0,n}=1$ and $\widetilde{c}_{1,n}=-c_{1,n+1}$ and the lemma follows.
\end{proof}

\begin{remark}
    For ${\sf w}(z)\equiv 1$, the value of $h_n$ specialises and all the odd-indexed coefficients vanish, that is, $c_{j,n}=\widetilde{c}_{j,n}=0$ for $j\geq 1$ odd.
\end{remark} 
Let us now define a RHP such that the 11 entry of its solution is the polynomial $\pi_n(z)$.
\begin{rhp}\label{rhp:Yn}
The Riemann-Hilbert Problem comprises of finding a $2\times 2$ matrix valued function $Y_n(z,\tau)$ with the following properties:
\begin{itemize}
    \item $Y_n(z,\tau)$ is analytic in $z\in \mathbb{T}\setminus ( \gamma \cup \{0\} )$.
    \item The following jump condition hold for $z\in \gamma$:
    \begin{align}\label{jump:Yn}
        Y_{n,+}(z,\tau) = Y_{n,-}(z,\tau) \left(\begin{array}{cc}
          1   &  {\sf w}(z)\\
          0   & 1
        \end{array} \right),
    \end{align}
    where, following the standard notation,
    $\pm$ indicate the piece-wise analytic functions to the left and right side respectively of $\gamma$ w.r.t its orientation, see Figure \ref{fig:1}.
    \item In the limit $z\to 0$:
    \begin{align}
        Y_n(z,\tau)=\left( \mathbb{1} + \mathcal{O}(z) \right) \left( \begin{array}{cc}
         z^{-n}    & 0 \\
           0  & z^{n-2}
        \end{array}\right).
    \end{align}
\end{itemize}
\end{rhp}


\begin{theorem}\label{thm:existence_uniqueness} Let $n\geq 3$, then RHP \ref{rhp:Yn} is uniquely solvable if and only if the determinant $D_n\neq 0$, in which case the solution is given by
\begin{align}
            Y_n(z,\tau) = \left(\begin{array}{cc}
             \pi_{n}(z) & \mathcal{C}(\pi_{n})(z)  \\
             \frac{2\pi i}{h_{n-1}}  \pi_{n-1}(z)    &  \frac{2\pi i}{h_{n-1}}  \mathcal{C}(\pi_{n-1})(z)
            \end{array} \right), \qquad (n\geq 3). \label{def:Yn}
\end{align}
\end{theorem}






\begin{proof}
To prove the theorem, we analyse the two rows of a solution to the RHP separately.

We start with the first row, $(Y_{11}(z),Y_{12}(z))$. The conditions imposed in the RHP translate to
\begin{itemize}
    \item $Y_{11}(z)$ and $Y_{12}(z)$ are analytic and doubly periodic on $\mathbb{T}\setminus ( \gamma \cup \{0\} )$,
    \item $Y_{11}(z)$ is analytic on $\gamma$ and $Y_{12}(z)$ satisfies the jump condition
    \begin{equation*}
        (Y_{12})_+(x)=(Y_{12})_-(x)+Y_{11}(x){\sf w}(x)\qquad (x\in \gamma).
    \end{equation*}
    \item asymptotic conditions
    \begin{equation*}
        Y_{11}(z)=z^{-n}(1+\mathcal{O}(z)),\quad Y_{12}(z)=\mathcal{O}(z^{n-1});\qquad z\rightarrow 0. 
    \end{equation*}
\end{itemize}
We refer to this as the row one RHP.

Note that, $Y_{11}(z)$ is necessarily a monic elliptic polynomial of degree $n$. In order to study $Y_{12}(z)$, we first consider the following function,
\begin{equation*}
    r(z):=Y_{12}(z)-\mathcal{C}(Y_{11})(z)+\zeta(z)\int_{\frac{\tau}{2}}^{\frac{\tau}{2}+1}Y_{11}(x){\sf {\sf w}}(x)dx.
\end{equation*}
Note that $r(z)$ has no jump on $\gamma$ and it is periodic with respect to $1$ and $\tau$, so it must be an elliptic function. Furthermore, by the asymptotic conditions,
\begin{equation*}
    r(z):=z^{-1}\int_{\frac{\tau}{2}}^{\frac{\tau}{2}+1}Y_{11}(x){\sf w}(x)dx+\mathcal{O}(1)\qquad (z\rightarrow 0).
\end{equation*}
Since $r(z)$ has no poles away from the lattice $\mathbb{Z}+\mathbb{Z}\cdot \tau$ and there exists no elliptic function of degree one, it follows that $r(z)$ must be a constant and
\begin{equation}\label{eq:ortho0}
    \int_{\frac{\tau}{2}}^{\frac{\tau}{2}+1}Y_{11}(x){\sf w}(x)dx=0.
\end{equation}

But, by the last equality, and the asymptotics of $Y_{11}(z)$ and $Y_{12}(z)$, we now have
\begin{equation*}
    r(z)=Y_{12}(z)-\mathcal{C}(Y_{11})(z)=\mathcal{O}(z^{n-1})
\end{equation*}
as $z\rightarrow 0$. Thus, $r(z)\equiv 0$ and it follows that
\begin{equation*}
    Y_{12}(z)=\mathcal{C}(Y_{11})(z).
\end{equation*}
In particular, by the asymptotics of $Y_{12}(z)$
\begin{equation*}
    \mathcal{C}(Y_{11})(z)=\mathcal{O}(z^{n-1})\quad (z\rightarrow 0).
\end{equation*}
In other words,
\begin{equation*}
     \int_{\frac{\tau}{2}}^{\frac{\tau}{2}+1}Y_{11}(z)\zeta^{(k)}(z){\sf w}(z)dz=0
\end{equation*}
for $1\leq k\leq n-2$, implying that
\begin{equation}\label{eq:y11ortho}
    \int_{\frac{\tau}{2}}^{\frac{\tau}{2}+1}Y_{11}(z)\mathcal{E}_m(z) {\sf w}(z)dz=0
    \end{equation}
for $2\leq m\leq n-1$. We have already seen that this equality also holds for $m=0$, see equation \eqref{eq:ortho0}.

Now, let us write  
\begin{equation*}
    Y_{11}(z)=c_0\mathcal{E}_0+c_2\mathcal{E}_2+c_3\mathcal{E}_3+\ldots+c_{n-1}\mathcal{E}_{n-1}+\mathcal{E}_n.
\end{equation*}
Then equation \eqref{eq:y11ortho} is equivalent to
\begin{equation}\label{eq:coefficients}
    S_n\begin{pmatrix} c_0\\
    c_2\\
    \vdots\\
    c_{n-1}\end{pmatrix}=\begin{pmatrix} -\mu_{n,0}\\
    -\mu_{n,2}\\
    \vdots\\
    -\mu_{n,n-1}\end{pmatrix},
\end{equation}
where $S_n$ is the matrix defined in equation \eqref{eq:Dn}.

All in all, we find that a solution $(Y_{11}(z),Y_{12}(z))$ to the row one RHP exists and is unique if and only if equation \eqref{eq:coefficients} has a unique solution, which in turn is true if and only if $D_n\neq 0$. In particular, in that case, $\pi_n(z)$ exists and
\begin{equation*}
    Y_{11}(z)=\pi_n(z),\quad Y_{12}(z)=\mathcal{C}(\pi_n)(z).
\end{equation*}

The second row of $Y$ is analysed in much the same way. We note that the second row $(Y_{21}(z),Y_{22}(z))$ must satisfy
\begin{itemize}
    \item $Y_{21}(z)$ and $Y_{22}(z)$ are analytic and doubly periodic on $\mathbb{T}\setminus ( \gamma \cup \{0\} )$,
    \item $Y_{21}(z)$ is analytic on $\gamma$ and $Y_{22}(z)$ satisfies the jump condition
    \begin{equation*}
        (Y_{22})_+(x)=(Y_{22})_-(x)+Y_{21}(x){\sf w}(x)\qquad (x\in \gamma).
    \end{equation*}
    \item asymptotic conditions
    \begin{equation*}
        Y_{21}(z)=\mathcal{O}(z^{-(n-1)})),\quad Y_{22}(z)=z^{n-2}(1+\mathcal{O}(z))\qquad (z\rightarrow 0). 
    \end{equation*}
\end{itemize}
We refer to this as the row two RHP.

It follows that $Y_{21}(z)$ is required to be an elliptic polynomial of degree less or equal to $n-1$, 
\begin{equation*}
    Y_{22}(z)=\mathcal{C}(Y_{21})(z),
\end{equation*}
and we find that
\begin{equation}\label{eq:y21ortho}
    \int_{\frac{\tau}{2}}^{\frac{\tau}{2}+1}Y_{21}(z)\mathcal{E}_m(z) {\sf w}(z)dz=\delta_{m,n-1}2\pi i,
    \end{equation}
for $m=0,2,3,\ldots, n-1$. Upon writing 
\begin{equation*}
    Y_{21}(z)=c_0\mathcal{E}_0+c_2\mathcal{E}_2+c_3\mathcal{E}_3+\ldots+c_{n-1}\mathcal{E}_{n-1},
\end{equation*}
this is equivalent to
\begin{equation}\label{eq:coefficients21}
    S_n\begin{pmatrix} c_0\\
    c_2\\
    \vdots\\
    c_{n-2}\\
    c_{n-1}\end{pmatrix}=\begin{pmatrix} 0\\
    0\\
    \vdots\\
    0\\
    2\pi i\end{pmatrix},
\end{equation}
where $S_n$ is the matrix defined in equation \eqref{eq:Dn}.

All in all, we find that a solution $(Y_{21}(z),Y_{22}(z))$ to row two RHP exists and is unique if and only if equation \eqref{eq:coefficients21} has a unique solution, which in turn is true if and only if $D_n\neq 0$. 

Furthermore, if $D_n\neq 0$, then 
\begin{equation*}
    \frac{2\pi i}{h_{n-1}}  \pi_{n-1}(z)=\frac{2\pi i}{D_n}\;\widetilde{\pi}_{n-1}(z),
\end{equation*}
is well-defined and necessarily
\begin{equation*}
    Y_{21}(z)=\frac{2\pi i}{h_{n-1}}  \pi_{n-1}(z),\quad Y_{22}(z)=\frac{2\pi i}{h_{n-1}}  \mathcal{C}(\pi_{n-1})(z).
\end{equation*}

We have now shown that both rows of a solution to the RHP exist if and only if $D_n\neq 0$, and, in such case, we proven that formula \eqref{rhp:Yn} holds. The theorem follows.
\end{proof}

The asymptotic expansions in Lemma \ref{Series lemma} now help us obtain the precise form of the determinant of the solution of the RHP: \ref{rhp:Yn} $Y_n$.
\begin{lemma}
    For ${\sf w}(z)\equiv 1$, the determinant of the solution takes the form
\begin{align}\label{eq:DetYn}
    \det Y_n(z,\tau) =
    \wp(z,\tau) + \alpha_n(\tau)=:  f_n(z),
\end{align}
where $\alpha_n(\tau)$ is given by
\begin{align}\label{def:alphabetan}
 \alpha_n := c_{2,n}+ \widetilde{c}_{2,n-1} - \beta_n, &&   \beta_n := \frac{h_n}{h_{n-1}}.
\end{align}
\end{lemma}
\begin{proof}
Setting ${\sf w}(z)\equiv 1$, let us begin by noting that $\det Y_n$
\begin{enumerate}
\item is a doubly-periodic function
\begin{align}\label{property1:detyn}
    \det Y_{n}(z+\tau) = \det Y_{n}(z+1) = \det Y_n(z),
\end{align}
\item and has no jump on $\gamma$,
\begin{align}\label{property2:detyn}
    \det Y_{n,+} = \det Y_{n,-}.
\end{align}
\end{enumerate}
With the asymptotic behaviour of $\pi_n(z)$ \eqref{asymp:pin}, and $\mathcal{C}(\pi_n)(z)$ \eqref{asymp:Cpinmixed} in the determinant of $Y_n$ \eqref{def:Yn}, we see that in the limit $z\to 0$
\begin{align}
    \det Y_n(z,\tau) = \frac{1}{z^2}+c_{2,n}+ \widetilde{c}_{2,n-1} - \frac{h_n}{h_{n-1}}+ \mathcal{O}(z^2).
\end{align}
Since $\det Y_n$ is an elliptic function due to properties \eqref{property1:detyn}-\eqref{property2:detyn}, it only has one pole, which is of order 2 at $z=0$, 
and therefore must be equal to the Weierstrass $\wp(z)$-function  plus a constant $\alpha_n$ defined in \eqref{def:Yn}, \eqref{def:alphabetan}.
\end{proof}
Let us now restrict to even polynomials and repeat the methods presented above.

\subsection{Hankel determinants and the even case}
Consider the vector space of even elliptic polynomials,
\begin{equation*}
    \mathcal{P}_{\text{even}}=\mathcal{P}_{\text{even}}(\tau)=\{\text{even elliptic functions with only a pole at $0$}\},
\end{equation*}
with corresponding basis
\begin{equation*}
    \{1,\wp(z),\wp(z)^2,\wp(z)^3,\ldots\}.
\end{equation*}
Let ${\sf w}$ be an even weight function on the interval $\tfrac{1}{2}\tau+[0,1]$, that is,
\begin{equation*}
    {\sf w}(\tfrac{1}{2}+\tfrac{1}{2}\tau+x)={\sf w}(\tfrac{1}{2}+\tfrac{1}{2}\tau-x),\qquad x\in \tfrac{1}{2}\tau+[0,1].
\end{equation*}
For any natural number $k\in\mathbb{N}$, we define the $k$th orthogonal polynomial $\pi_{2k}(z)$ with respect to ${\sf w}$, if it exists,  by the conditions
\begin{align*}
    &\int_{\frac{1}{2}\tau+0}^{\frac{1}{2}\tau+1}\wp(z)^m\pi_{2k}(z){\sf w}(z)dz=0,\quad (0\leq m<2k),\\
    &\pi_{2k}(z)=\wp(z)^k(1+\mathcal{O}(\wp(z)^{-1}))\quad (z\rightarrow 0).
\end{align*}
For $i+j=k\in\mathbb{N}$, 
\begin{align}\label{2kth_moment}
    \nu_{2i,2j}\equiv  \nu_{2(i+j)} = \nu_{2k}:=\int_{\frac{1}{2}\tau+0}^{\frac{1}{2}\tau+1}\wp(z)^{k} {\sf w}(z)dz.
\end{align}
Then, $\pi_{2k}(z)$ exists if and only if the Hankel determinant of moments
\begin{equation}\label{eq:delta2k}
 \Delta_{2k}:=\begin{vmatrix} 
 \nu_{0} & \nu_{2} & \ldots & \nu_{2k-2}\\
 \nu_{2} & \nu_{4} &\ldots  & \nu_{2k}\\
 \vdots & \vdots  &\ddots & \vdots\\
 \nu_{2k-2} & \nu_{2k} & \ldots & \nu_{4k-2}
 \end{vmatrix},
 \end{equation}
 is nonzero, with $\Delta_0=1$. In turn, $\pi_{2k}$ is explicitly given by
 \begin{equation*}
\pi_{2k}(z)=\frac{1}{\Delta_{2k}}\begin{vmatrix} 
 \nu_{0} & \nu_{2} & \ldots & \nu_{2k}\\
 \nu_{2} & \nu_{4} &\ldots  & \nu_{2k+2}\\
 \vdots & \vdots  &\ddots & \vdots\\
 \nu_{2k-2} & \nu_{2k} & \ldots & \nu_{4k-2}\\
 1 & \wp(z) &\ldots & \wp(z)^k
 \end{vmatrix},
 \end{equation*}
further implying that
\begin{equation}\label{eq:normdelta}
h_{2k}:=\int_{\frac{1}{2}\tau+0}^{\frac{1}{2}\tau+1}\pi_{2k}(z)^2{\sf w}(z)dz=\int_{\frac{1}{2}\tau+0}^{\frac{1}{2}\tau+1}\wp(z)^k\pi_{2k}(z){\sf w}(z)dz=\frac{\Delta_{2k+2}}{\Delta_{2k}}.
\end{equation}
Another direct consequence is that the even polynomial can now be expanded as
\begin{equation}\label{eq:OPexpansion}
    \pi_{2k}(z)=\wp(z)^k-\frac{\Gamma_{2k}}{\Delta_{2k}}\wp(z)^{k-1}+\frac{\Lambda_{2k}}{\Delta_{2k}}\wp(z)^{k-2}+\mathcal{O}(\wp(z)^{k-3}),
\end{equation}
where
\begin{align*}
 \Gamma_{2k}&:=\begin{vmatrix} 
 \nu_{0} & \nu_{2} & \ldots & \nu_{2k-4} & \nu_{2k}\\
 \nu_{2} & \nu_{4} &\ldots  & \nu_{2k-2} & \nu_{2k+2}\\\
 \vdots & \vdots  &\ddots & \vdots & \vdots\\
 \nu_{2k-2} & \nu_{2k} & \ldots & \nu_{4k-6} & \nu_{4k-2}
 \end{vmatrix},&&
 \Lambda_{2k}&:=\begin{vmatrix} 
 \nu_{0} & \nu_{2} & \ldots & \nu_{2k-6} & \nu_{2k-2} & \nu_{2k}\\
 \nu_{2} & \nu_{4} &\ldots  & \nu_{2k-4} & \nu_{2k} & \nu_{2k+2}\\
 \vdots & \vdots  &\ddots & \vdots & \vdots & \vdots\\
 \nu_{2k-2} & \nu_{2k} & \ldots & \nu_{4k-8} & \nu_{4k-4} & \nu_{4k-2}
 \end{vmatrix},
 \end{align*}
for $k\geq 1$ with $\Gamma_2=\mu_2$ and $\Gamma_0=\Lambda_2=\Lambda_0=0$. 

\begin{proposition}
    For ${\sf w}(z)\equiv 1$, the moments satisfy the following recursion
    \begin{align}\label{nu2k_rec}
          (8k+12)\nu_{2k+4}=(2k+1)g_2\nu_{2k}+2kg_3\nu_{2k-2}.
    \end{align}
\end{proposition}
\begin{proof}
With the prescribed weight function, the moment \eqref{2kth_moment} reads
    \begin{equation}\label{eq:even_xy1}
   \nu_{2k}=\int_{\frac{1}{2}\tau+0}^{\frac{1}{2}\tau+1}\wp(z)^kdz=\oint_{\text{cycle}}\frac{x^k}{y}dx,
\end{equation}
where $y^2=4x^3-g_2 x-g_3$, as can be seen by a change of variables $x=\wp(z)$. Note that expression
\begin{equation*}
    \frac{d}{dx}x^k y=kx^{k-1}y+\frac{1}{2y}(12 x^{k+2}-g_2 x^k),
\end{equation*}
from which we obtain the following identity
\begin{equation}\label{eq:identity}
    k\oint_{cycle} x^{k-1}y\,dx=\tfrac{g_2}{2} \nu_{2k}-6\nu_{2k+4},
\end{equation}
for $k\geq 1$. Therefore,
\begin{align}
4\nu_{2k+4}=&4\oint_{\text{cycle}}\frac{4x^{k+2}}{y}dx,\\
    =&\oint_{\text{cycle}}\frac{x^{k-1} y^2+g_2 x^{k}+g_3x^{k-1}}{y}dx\\
    =&\oint_{\text{cycle}}x^{k-1} ydx+g_2 \nu_{2k}+g_3\nu_{2k-2}\\
    =&\frac{1}{k}(\tfrac{g_2}{2} \nu_{2k}-6\nu_{2k+4})+g_2 \nu_{2k}+g_3\nu_{2k-2},
\end{align}
where we used identity \eqref{eq:identity} in the last equality. This gives the recursion in the proposition.
\end{proof}
The first values are
\begin{align*}
    \nu_0=1,&&
    \nu_1=-2\eta_1,&&
    \nu_2=\frac{1}{12}g_2, && \nu_3=\frac{1}{10}(g_3-3g_2\eta_1),
\end{align*}
where $\eta_1=\zeta(\frac{1}{2})$ is the first period of the second kind. Correspondingly, the first few even Hankel determinants, defined in equation \eqref{eq:delta2k}, are given by
\begin{align*}
    &\Delta_0=\Delta_2=1,\\
    &\Delta_4=\frac{1}{12}(g_2-48 \eta_1^2),\\
    &\Delta_6=\frac{1}{37800}(25 g_2^3-378 g_3^2+108 g_2 g_3\eta_1-1872 g_2^2\eta_1^2+43200 g_3\eta_1^3).
\end{align*}

The Hankel determinants of moments are functions of the modular parameter $\tau\in\mathbb{H}$, and they can be written as explicit functions in $\tau$ using the equations for $g_2,g_3$ and $\eta_1$, as functions of $\tau$, in Appendix \ref{Appendix:EllFunc}. They satisfy the following symmetries,
\begin{equation*}
    \Delta_{2k}(\tau)^\diamond=\Delta_{2k}(\tau^\diamond),\quad \Delta_{2k}(\tau+2)=\Delta_{2k}(\tau)\quad (\tau\in\mathbb{H}),
\end{equation*}
where $(\alpha+\beta i)^\diamond=-\alpha+\beta i$ for $\alpha,\beta\in\mathbb{R}$.

\begin{figure}[h]\captionsetup[subfigure]{labelformat=empty}
	\centering
\begin{subfigure}[b]{0.46\textwidth}
\centering
 \includegraphics[width=\textwidth]{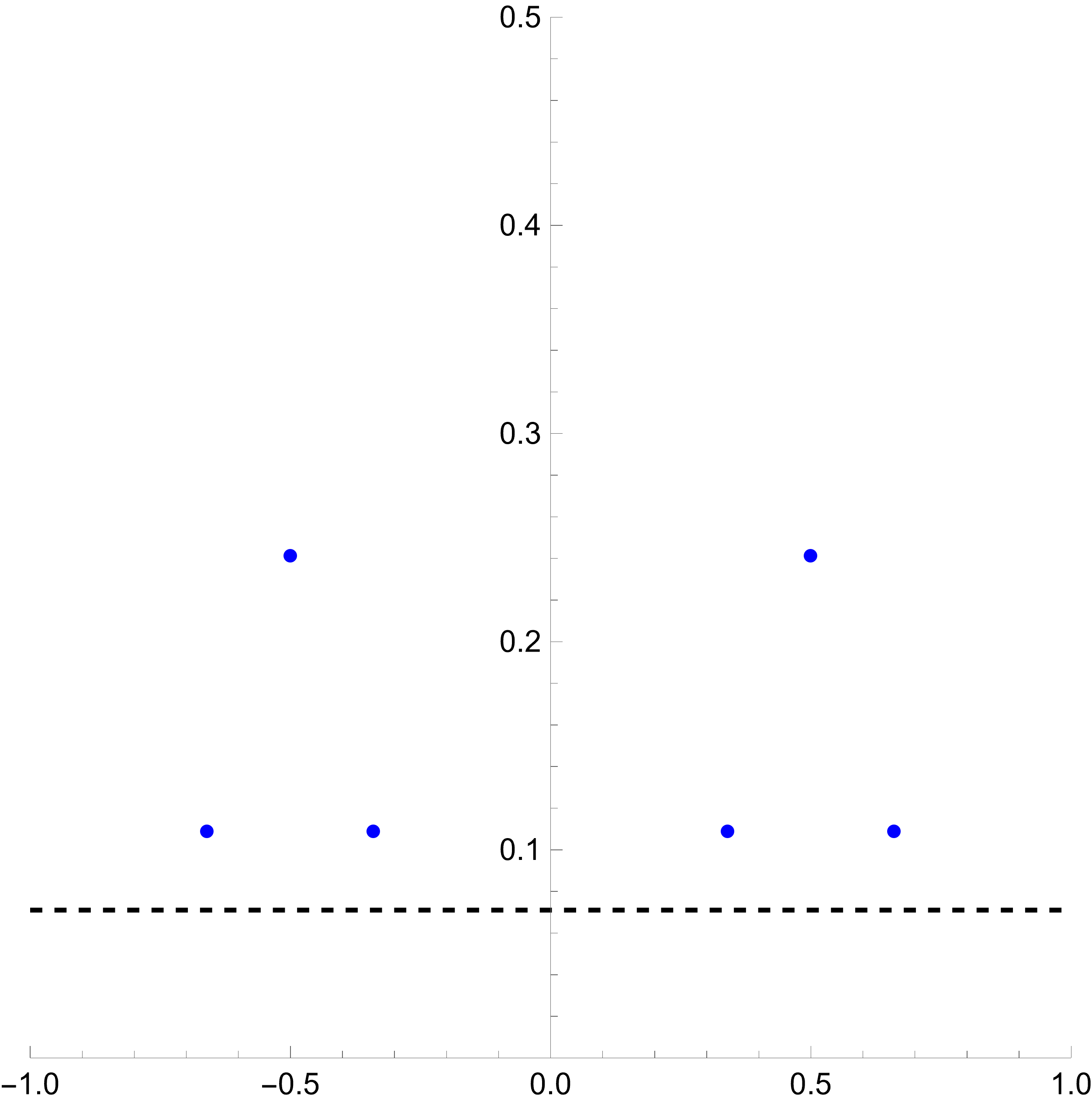}
 \caption{$k=2$}
 \end{subfigure}\hfill
 \begin{subfigure}[b]{0.46\textwidth} 
 \centering
 \includegraphics[width=\textwidth]{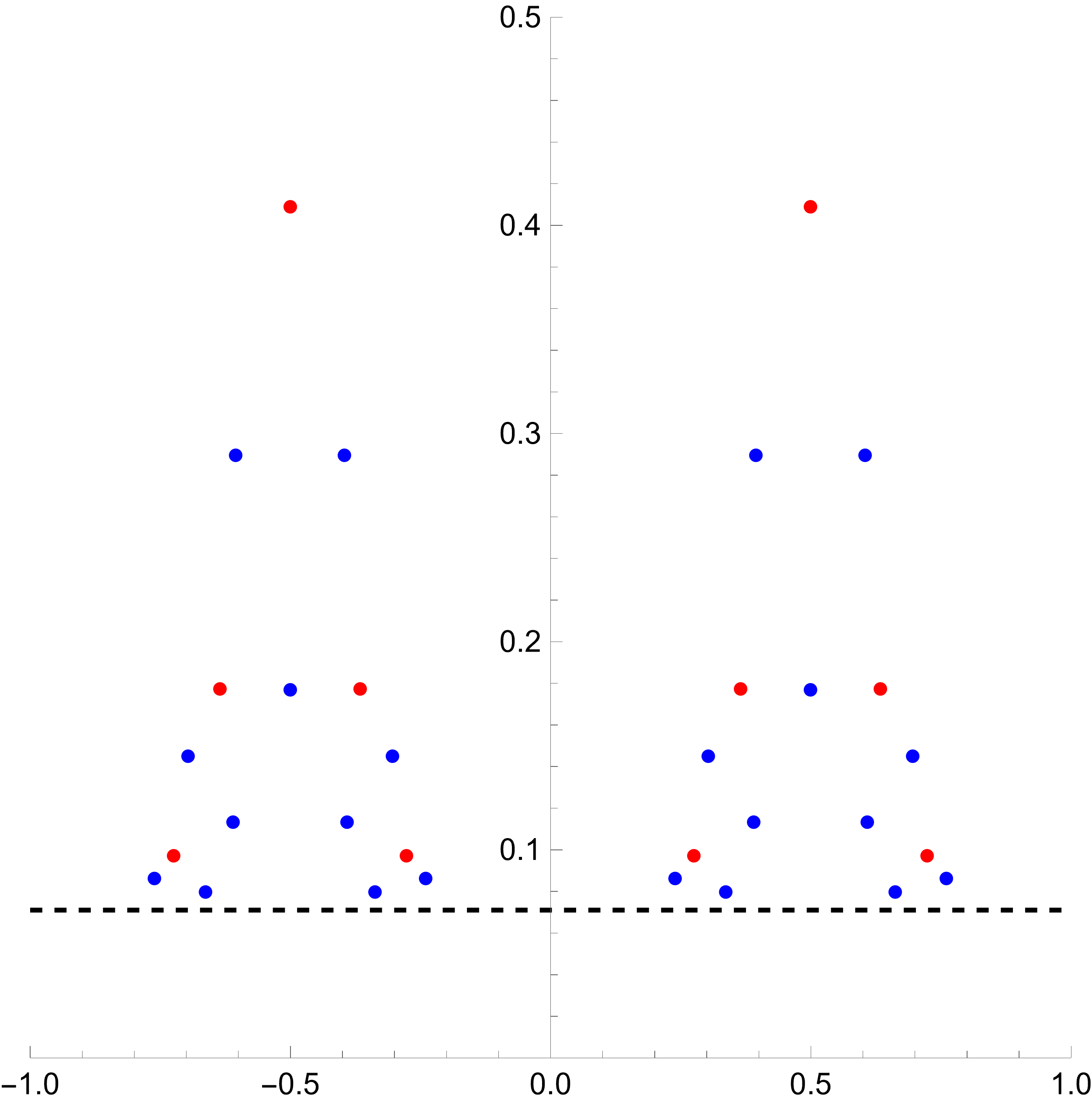} 
 \caption{$k=3$}
 \end{subfigure}\vspace{1cm}
  \begin{subfigure}[b]{0.46\textwidth}
  \centering
 \includegraphics[width=\textwidth]{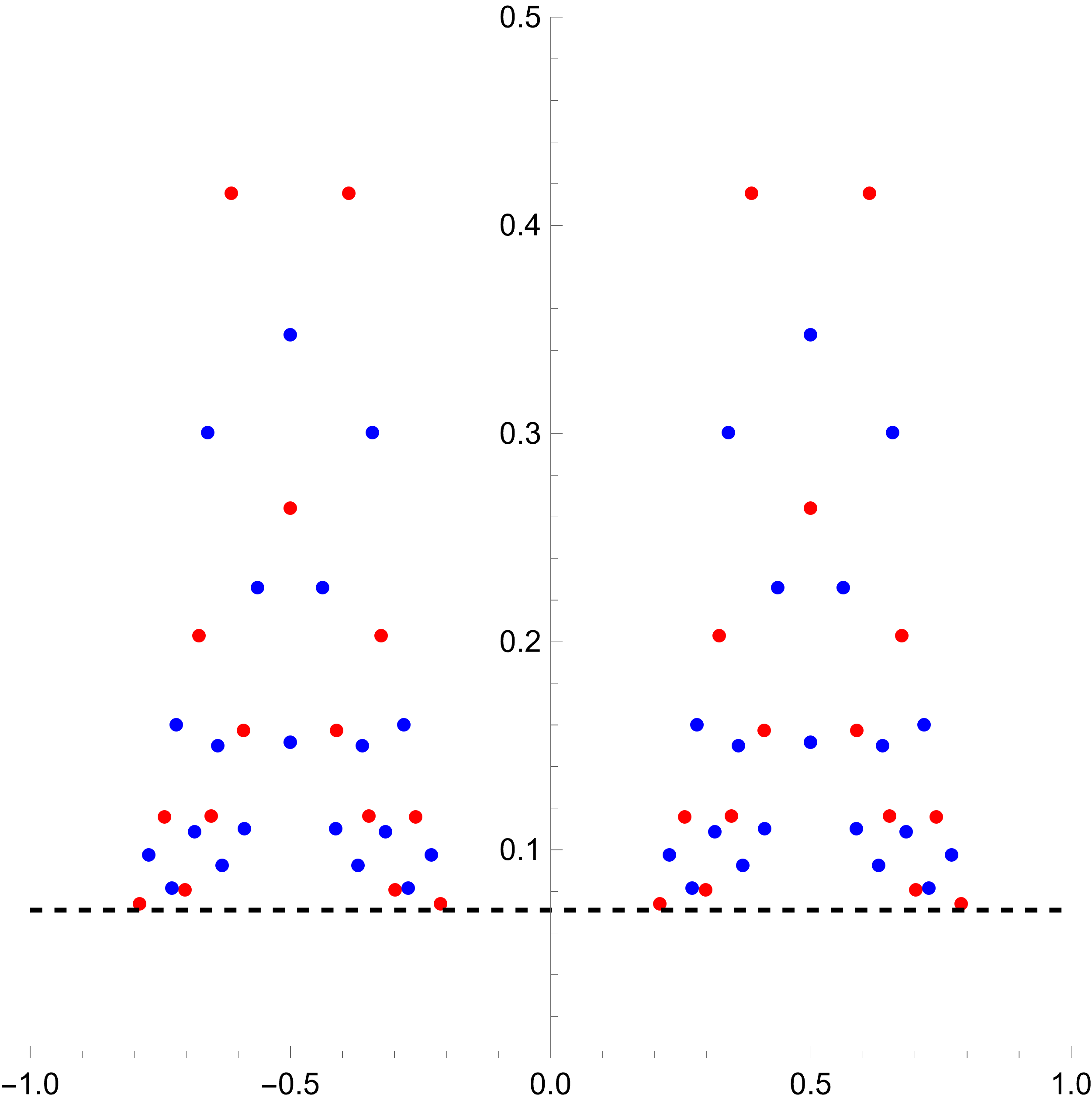}
 \caption{$k=4$}
 \end{subfigure}
 \hfill \begin{subfigure}[b]{0.46\textwidth}
 \centering
 \includegraphics[width=\textwidth]{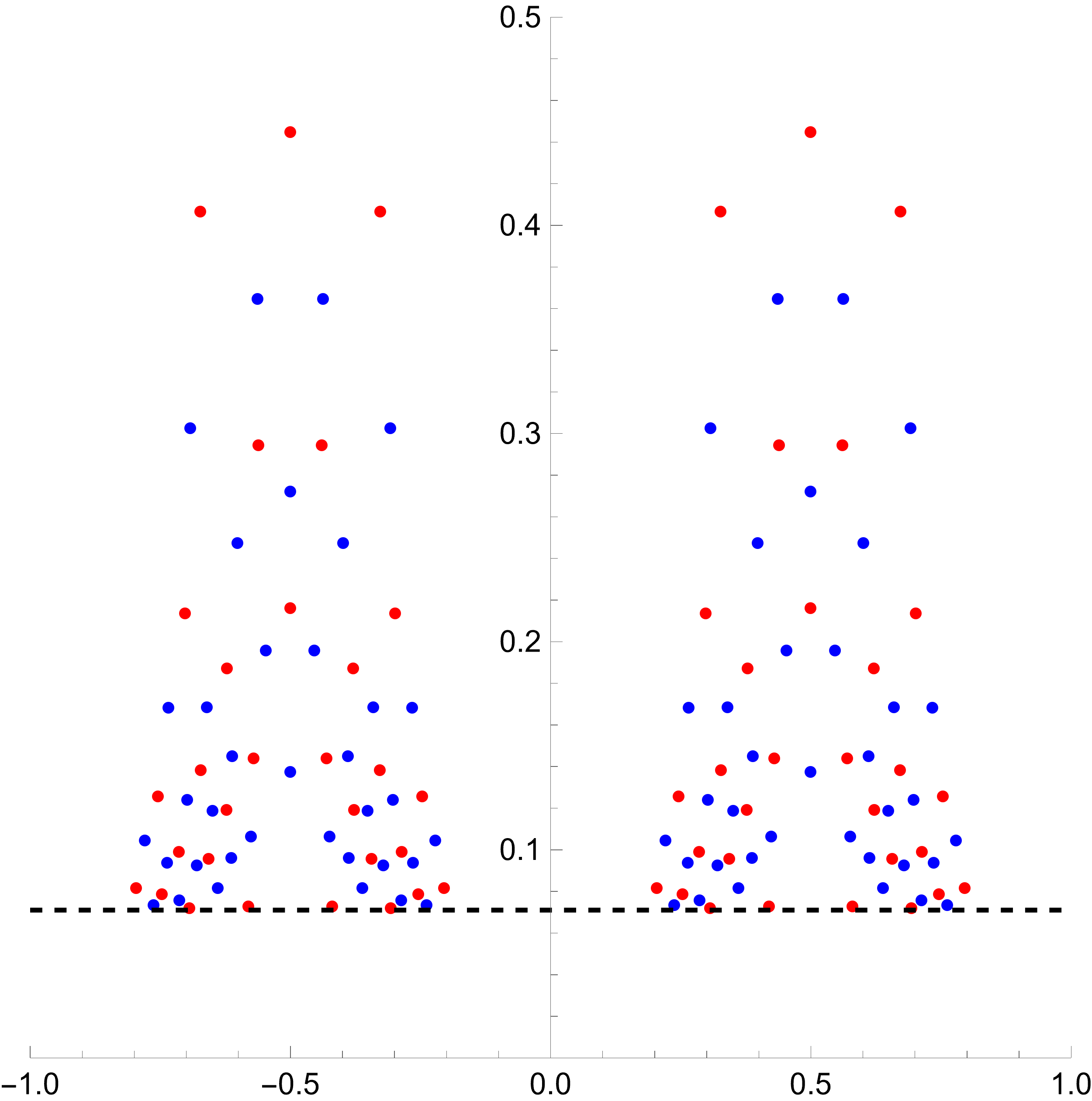}
  \caption{$k=5$}
 \end{subfigure}
\caption{In these plots the zero distributions of $\Delta_{2k-1}(\tau)$ and $\Delta_{2k}(\tau)$ are displayed in respectively red and blue, for $k=2,3,4,5$, in the domain defined by $-1\leq \Re\;\tau\leq 1$ and $\Im \; \tau\geq \tfrac{1}{\pi} \log(\tfrac{5}{4})\approx 0.071$. Furthermore, in each plot the dashed line is the line $\Im \; \tau=\tfrac{1}{\pi} \log(\tfrac{5}{4})$.
} \label{fig:zero_distributions}
\end{figure}

In Figure \ref{fig:zero_distributions}, the zero distribution of $\Delta_{2k}(\tau)$ is displayed in blue, for $k=2,4,6,8$, in the upper-half plane cut off by
\begin{equation*}
 \Im \; \tau\geq \tfrac{1}{\pi} \log(\tfrac{5}{4}).   
\end{equation*}
The reason for this cut, is that the numerics become unstable near the real line. In particular, even though only finitely many zeros are shown in the plots, there might in fact be an infinite number of zeros accumulating at points $\tau\in\mathbb{Z}$, for the even Hankel determinants $\Delta_{2k}$, $k\geq 2$.

\begin{remark}
A similar computation follows for the odd case starting from \eqref{def:moment}:
\begin{align}
    \nu_{2i+3, 2j+3} = \int_{\frac{\tau}{2}}^{1+\frac{\tau}{2}} (\wp'(z))^2 \wp(z)^{i+j} dz = \oint_{\text{cycle}} y x^k dx
\end{align}
with the change of variables as in \eqref{eq:even_xy1}. After a direct manipulation, we get
\begin{align}
    \nu_{2k+3} = 4 \nu_{2k+6} - g_2 \nu_{2k+4}-g_3 \nu_{2k}.
\end{align}
\end{remark}
The determinant of the odd moments is
\begin{equation*}
 \Delta_{2k+3}:=\begin{vmatrix} 
 \nu_{3} & \nu_{5} & \ldots & \nu_{2k+1}\\
 \nu_{5} & \nu_{7} &\ldots  & \nu_{2k+3}\\
 \vdots & \vdots  &\ddots & \vdots\\
 \nu_{2k+1} & \nu_{2k+3} & \ldots & \nu_{4k-1}
 \end{vmatrix}.
\end{equation*}
We have plotted the zero distributions of the first couple of odd Hankel determinants in red in Figure \ref{fig:zero_distributions}.

Note that $\Delta_{2k}$ and $\Delta_{2k+1}$ are Hankel determinants, for $k\geq 0$. However, the determinant $D_n$ is generally not Hankel, but has a checkerboard pattern of even and odd moments leading to the factorisation
     \begin{align*}
         D_n = \Delta_{n} \Delta_{n+1}, && n\geq 2.
     \end{align*}

\subsection{Riemann-Hilbert problem for the even case}
If the weight function is even on $\gamma$, i.e. equation \eqref{eq:even_weight} holds, the sequence of orthogonal polynomials splits into two sequences of even and odd EOPs  and we can define a RHP corresponding to the even polynomials $\pi_{2k}$, $k\geq 0$. The analysis for the odd case mirrors the even one and so we restrict our analysis to the even case. 

Analogous to RHP \ref{rhp:Yn}, the even polynomials $\pi_{2k}$ appear as the 1,1 elements of the solution 
 \begin{align}
           Y_{2k}(z,\tau) &= \left(\begin{array}{cc}
    \pi_{2k}(z) & \mathcal{C}(\pi_{2k}{\sf w})(z)  \\
    \frac{2\pi i}{h_{2k-2}}  \pi_{2k-2}(z)    &  \frac{2\pi i}{h_{2k-2}}  \mathcal{C}(\pi_{2k-2}{\sf w})(z)
    \end{array} \right), && k\geq 1, \label{rhpsol:even}
        \end{align}
        to the following RHP. 
\begin{rhp}\label{rhp:even}
    \begin{itemize}
        \item The function $Y_{2k}(z,\tau)$ is piece-wise analytic on $\mathbb{T}\backslash\left(\gamma\cup \lbrace 0\rbrace\right)$,
        \item for $z\in \gamma$, the following jump condition holds
        \begin{align}
             Y_{2k,+}(z,\tau) = Y_{2k,-}(z,\tau) \left( \begin{array}{cc}
                 1 & {\sf w}(z) \\
                 0 & 1
             \end{array}\right),
        \end{align}
        \item and in the limit $z\to 0$,
        \begin{align}\label{asymp:Y2k}
             Y_{2k}(z,\tau) &= \left(\mathbb{1}+\mathcal{O}(z) \right) \left(\begin{array}{cc}
        z^{-2k} & 0 \\
       0  & z^{2k-3}
    \end{array} \right).
        \end{align}
    \end{itemize}
\end{rhp}
Analogously to Theorem \ref{thm:existence_uniqueness},
RHP is \ref{rhp:even} is solvable if and only if the Hankel determinant of even moments $\Delta_{2k}(\tau)\neq 0$.

Mimicking the analysis for the general polynomials, we begin with the asymptotic behaviour of the even polynomials and their Cauchy transform.
\begin{lemma}\label{lemma:evenpol}
    In the limit $z\to 0$, the polynomials 
    \begin{align}
        \pi_{2k}(z)= \sum_{i=0}^{\infty} \frac{a_{i,2k}}{z^{2(k-i)}}, \qquad a_{0,2k} =1, 
    \end{align}
    and their Cauchy transform 
    \begin{align}\label{asymp:Cpin}
        \mathcal{C}(\pi_{2k})(z) = \frac{h_{2k}}{2\pi i} \sum_{i=0}^{\infty} \widetilde{a}_{i,2k} z^{2(k+i)-1}, \qquad \widetilde{a}_{0,2k} =1, && \widetilde{a}_{1,2k} =-a_{1, 2k+2}.
    \end{align}
\end{lemma}
\begin{proof}
    The proof is the same as for Lemma \ref{Series lemma}.
\end{proof}
\begin{lemma}\label{lemma:detY2k}
The determinant of $Y_{2k}$ is
    \begin{align}
    \det Y_{2k}(z,\tau) = -\frac{ \wp'(z)}{ 2}. \label{eq:DetY2k}
\end{align}
\end{lemma}
\begin{proof}
Consider the determinant $\det Y_{2k}(z)$ of the solution to RHP \ref{rhp:even}. Note that $\det Y_{2k}(z)$ has a trivial jump along $\gamma$ and thus extends to an elliptic function with only a pole at $z=0$. Furthermore, from the asymptotic behaviour of $Y_{2k}(z)$, it follows that
\begin{equation*}
    \det Y_{2k}(z)=z^{-3}(1+\mathcal{O}(z))\qquad (z\rightarrow 0).
\end{equation*}
Finally, note that $Y_{2k}(-z)\sigma_3$ also solves RHP \ref{rhp:even}, thus $Y_{2k}(z)=Y_{2k}(-z)\sigma_3$, and $\det Y_{2k}(-z)=-\det Y_{2k}(z)$. We have now shown that $\det Y_{2k}(z)$ is an odd, monic, elliptic polynomial of degree $3$. There exists only one such polynomial, $-\frac{1}{2}\wp'(z)$. The lemma follows.
\end{proof}

\section{Linear problems for the general case and recurrence relations}\label{Sec:Mix_Lin}

In this section we show that the solution of the RHP \ref{rhp:Yn} for the general polynomials $\pi_n$ satisfies a system of linear differential and difference equations, and their compatibility condition leads to the recurrence relation for the polynomials and the recursion for their coefficients. Henceforth, we use the notation
\begin{align*}
    ' = \frac{d}{dz}, && \dot{} = \frac{d}{d\tau}
\end{align*}


\begin{proposition}\label{prop:difference_ls_Yn}
The solution of the RHP: \ref{rhp:Yn} solves the linear difference equation
\begin{align}\label{linsys:n+1}
    Y_{n+1} = R_n Y_n, && R_n = \frac{1}{f_n} \left(\begin{array}{cc}
     -\wp'(z)/2   & -\frac{h_{n}}{2\pi i } f_{n+1}  \\
      \frac{2\pi i}{h_n} f_n  & 0
    \end{array} \right), && n\geq 3,
\end{align}
where $f_n$ is the determinant of $Y_n$ \eqref{eq:DetYn}.
\end{proposition}
\begin{proof}
We begin with the $(n+1)^{th}$ solution \eqref{def:Yn} 
\begin{align}
    Y_{n+1}(z,\tau) = \left(\begin{array}{cc}
             \pi_{n+1}(z) & \mathcal{C}(\pi_{n+1})(z)  \\
             \frac{2\pi i}{h_{n}}  \pi_{n}(z)    &  \frac{2\pi i}{h_{n}}  \mathcal{C}(\pi_{n})(z)
            \end{array} \right),
\end{align}
and observe that
\begin{align}
 Y_{n+1} Y_{n}^{-1} \det Y_{n} &= \left(\begin{array}{cc}
             \pi_{n+1}(z) & \mathcal{C}(\pi_{n+1})(z)  \\
             \frac{2\pi i}{h_{n}}  \pi_{n}(z)    &  \frac{2\pi i}{h_{n}}  \mathcal{C}(\pi_{n})(z)
            \end{array} \right)\left(\begin{array}{cc}
            \frac{2\pi i}{h_{n-1}}  \mathcal{C}(\pi_{n-1})(z) & -\mathcal{C}(\pi_{n})(z)  \\
             -\frac{2\pi i}{h_{n-1}}  \pi_{n-1}(z)    &  \pi_n(z)
            \end{array} \right)\nonumber\\
            &=\left(\begin{array}{cc}
    \frac{2 \pi i }{h_{n-1}} \left(\pi_{n+1}  \mathcal{C}(\pi_{n-1})- \pi_{n-1} \mathcal{C}(\pi_{n+1}) \right)   & -\frac{h_{n}}{2 \pi i}\det(Y_{n+1})  \\
      \frac{2\pi i}{h_n} \det(Y_n)   & 0
    \end{array} \right).\label{comp:Yn+1Yninv}
\end{align}
The 11 element in the above expression is determined as follows. In the limit $z\to 0$,
\begin{align}
    \frac{2\pi i}{h_{n-1}} \left(\pi_{n+1}(z) \mathcal{C}(\pi_{n-1})(z)-\pi_{n-1}(z) \mathcal{C}(\pi_{n+1})(z)\right) = \frac{1}{z^3} \left( 1+\mathcal{O}(z) \right),
\end{align}
and following from the periodicity properties of $Y_n$, $\det Y_n$, the LHS of \eqref{comp:Yn+1Yninv} is a matrix valued elliptic function. Therefore,
\begin{align}
 \left( Y_{n+1} Y_{n}^{-1} \det Y_{n} \right)_{11}= -\frac{\wp'(z)}{2}.
\end{align}
The above expression along with \eqref{eq:DetYn} gives \eqref{linsys:n+1}.
\end{proof}

\begin{proposition}\label{prop:differential_ls_Yn}
    The solution of the RHP: \ref{rhp:Yn}  $Y_n$ solves the following linear differential equation for $n\geq 3$:
    \begin{align}\label{linsys:dz}
    Y_n' = L_n Y_n, && L_n = \frac{1}{f_n}\left(\begin{array}{cc}
        n \wp'(z)/2 & \frac{h_n}{2\pi i } \left(  (n-1) f_n + n f_{n+1}\right)  \\
        \frac{2\pi i}{h_{n-1}} \left( (2-n) f_{n-1} + (1-n) f_n \right) & (2-n) \wp'(z)/2
    \end{array} \right).
\end{align}
\end{proposition}
\begin{proof}
    With the derivative\footnote{We drop the $z$, $\tau$ dependence in favour of brevity.} of $Y_n$ in \eqref{def:Yn}, we get
\begin{align}\label{comp:YnpYninv}
     Y_{n}' Y_{n}^{-1} \det Y_{n} &= \left(\begin{array}{cc}
             \pi_{n}'(z) & \partial_z\mathcal{C}(\pi_{n})(z)  \\
             \frac{2\pi i}{h_{n}}  \pi_{n}'(z)    &  \frac{2\pi i}{h_{n}}  \partial_z\mathcal{C}(\pi_{n})(z)
            \end{array} \right)\left(\begin{array}{cc}
            \frac{2\pi i}{h_{n-1}}  \mathcal{C}(\pi_{n-1})(z) & -\mathcal{C}(\pi_{n})(z)  \\
             -\frac{2\pi i}{h_{n-1}}  \pi_{n-1}(z)    &  \pi_n(z)
            \end{array} \right).
\end{align}
Recalling the asymptotic behaviour of the following entities near $z\to 0$ in Lemma \ref{Series lemma} now paying attention to the fact that with our weight the polynomials split into odd and even parts,
\begin{align}
    \pi_n (z) &= \frac{1}{z^n} + \frac{c_{2,n}}{z^{n-2}}+\dots, \\
    \mathcal{C}(\pi_n) (z) &= \frac{h_n}{2\pi i} \left(z^{n-1}+\widetilde{c}_{2,n} z^{n+1}+\dots  \right),
\end{align}
the asymptotic behaviour of the 11 element of \eqref{comp:YnpYninv} is
\begin{align}
     \left( Y_{n}' Y_{n}^{-1} \det Y_{n}\right)_{11} &= -\frac{n}{z^3} -\frac{1}{z}\left( n \widetilde{c}_{2, n-1} + (n-2) c_{2, n} + (n-1) \beta_n \right)+\mathcal{O}(z),
\end{align}
and because the LHS is an elliptic function,
\begin{align}
     \left( Y_{n}' Y_{n}^{-1} \det Y_{n}\right)_{11} = n \frac{\wp'(z)}{2},
\end{align}
with the constraint
\begin{align}\label{eq:betan_c}
    n \widetilde{c}_{2, n-1} + (n-2) c_{2, n} + (n-1) \beta_n =0.
\end{align}
Also note that
\begin{align}\label{detcons1}
    c_{2,n} + \widetilde{c}_{2,n-2}=0.
\end{align}
The 12 element in the limit $z\to 0$ is
\begin{align}
 \left( Y_{n}' Y_{n}^{-1} \det Y_{n}\right)_{12}=    \frac{h_n}{2\pi i} \left(\frac{n}{z^2} + \frac{(n-1)}{z^2} + (2n+1) \widetilde{c}_{2,n} + (2n-3) c_{2,n} + \mathcal{O}(z^2)\right),
\end{align}
and we can check that 
\begin{align}
    (n-1) f_n + n f_{n+1} &= (n-1) \left(\wp(z) + c_{2,n} + \widetilde{c}_{2,n-1} - \beta_n \right) + n \left(\wp(z) + c_{2,n+1} + \widetilde{c}_{2,n} - \beta_{n+1} \right) \nonumber \\
    & \mathop{=}^{\eqref{eq:betan_c}}\frac{(n-1)}{z^2} + \frac{n}{z^2} + (n-1) c_{2,n} + (n-1) \widetilde{c}_{2,n-1} +n \widetilde{c}_{2,n-1} + (n-2)c_{2,n} \nonumber\\
    &+ n c_{2,n+1} + n \widetilde{c}_{2,n} +(n+1) \widetilde{c}_{2,n} + (n-1) c_{2,n+1} \nonumber \\
    &\mathop{=}^{\eqref{detcons1}}\frac{(n-1)}{z^2} + \frac{n}{z^2} + (2n-3) c_{2,n}  +(2n+1) \widetilde{c}_{2,n}.
\end{align}
Therefore,
\begin{align}
     \left( Y_{n}' Y_{n}^{-1} \det Y_{n}\right)_{12} = \frac{h_n}{2\pi i } \left(  (n-1) f_n + n f_{n+1}\right).
\end{align}
Similarly, the 21 element
\begin{align}
     \left( Y_{n}' Y_{n}^{-1} \det Y_{n}\right)_{21}=\frac{2\pi i}{h_{n-1}} \left(\frac{1-n}{z^2} + \frac{2-n}{z^2} + (5-2n) c_{2,n-1} + (1-2n) \widetilde{c}_{2, n-1} +\mathcal{O}(z^2)\right),
\end{align}
and as before, we can see that 
\begin{align}
    (2-n) f_{n-1} + (1-n) f_n &= (2-n) \left( \wp(z) + c_{2,n-1} + \widetilde{c}_{2,n-2} - \beta_{n-1} \right) \nonumber \\
    & + (1-n)\left( \wp(z) + c_{2,n} + \widetilde{c}_{2,n-1} - \beta_{n} \right) \nonumber \\
    & = \frac{(2-n)}{z^2} + \frac{(1-n)}{z^2} + (2-n) c_{2,n-1} + (2-n)\widetilde{c}_{2,n-2} - (2-n)\beta_{n-1} \nonumber \\
    & + (1-n)c_{2,n} + (1-n)\widetilde{c}_{2,n-1} - (1-n)\beta_{n} \nonumber \\
    & \mathop{=}^{\eqref{eq:betan_c}} \frac{(2-n)}{z^2} + \frac{(1-n)}{z^2} + (2-n) c_{2,n-1} + (2-n)\widetilde{c}_{2,n-2} + (1-n) \widetilde{c}_{2,n-2} \nonumber \\
    &+ (3-n) c_{2, n-1}  + (1-n)c_{2,n} + (1-n)\widetilde{c}_{2,n-1} -n \widetilde{c}_{2, n-1} + (2-n) c_{2,n} \nonumber \\
    & \mathop{=}^{\eqref{detcons1}} \frac{(2-n)}{z^2} + \frac{(1-n)}{z^2} + (5-2n) c_{2,n-1} + (1-2n) \widetilde{c}_{2,n-1}.
\end{align}
Therefore,
\begin{align}
     \left( Y_{n}' Y_{n}^{-1} \det Y_{n}\right)_{21}= \frac{2\pi i}{h_{n-1}} \left( (2-n) f_{n-1} + (1-n) f_n \right).
\end{align}
Finally, we see that 
\begin{align}
    \left( Y_{n}' Y_{n}^{-1} \det Y_{n}\right)_{22} = (2-n) \frac{\wp'(z)}{2},
\end{align}
which finishes the proof of the proposition.
\end{proof}
An immediate consequence of the above proposition is that the 11 entry of the linear equation \eqref{linsys:dz} gives an ODE for the polynomials.
\begin{corollary} The polynomials $\pi_n$ solve the following second order differential equation
  \begin{align*}
         \pi_n''(z)
         &=  \left(\frac{ \wp'(z)}{f_n}  + n  \left( \frac{f_{n+1}}{f_{n}}\right)'   \left( (n-1)  + n \frac{f_{n+1}}{f_{n}} \right)^{-1} \right)\pi_{n}' \\
         &+\left(\left(\frac{n \wp'(z)}{2f_n} \right)'-n  \left( \frac{f_{n+1}}{f_{n}}\right)'   \frac{n \wp'(z)}{2\left(  (n-1) f_n + n f_{n+1}\right)}- \det L_n\right)\pi_n. 
    \end{align*}
\end{corollary}


\subsection{Recurrence relations}

With the above linear system, as first step, we obtain the recurrence relation for the elliptic-polynomials.
We find
\begin{equation}
    \pi_{n+1} = -\frac{\wp'\pi_{n}}{2f_n}-\frac{\beta_nf_{n+1}}{f_n}\pi_{n-1},\label{scal lax}
\end{equation}
which follows from the first entry of \eqref{linsys:n+1}. Note that $f_n$ as defined in \eqref{eq:DetYn} depends on $z$. Iterating the above equation once again yields the following theorem.
\begin{theorem}\label{thm:rec_pol_gen}
The polynomials $\{\pi_n\}_{n=4}^\infty$, satisfy the relation
    \begin{equation}\label{eq:rec_pin2pin}
        \pi_{n+2} = \left(\wp - B_n \right)\pi_n - \beta_n\beta_{n-1}\pi_{n-2},
\end{equation}
where,
\[ B_n = \beta_{n+1}+\beta_n + \alpha_n + \alpha_{n+1}.\] 
Furthermore,
\begin{equation}\label{comrec}
    -\alpha_n^3 +\frac{g_2}{4}\alpha_n -\frac{g_3}{4} = \beta_n(\alpha_{n-1}-\alpha_n)(\alpha_{n+1}-\alpha_n).
\end{equation}
\end{theorem}

\begin{proof}
Iterating \eqref{scal lax} by $n\pm1$ we find that
\begin{equation}\label{eq 6}
    \pi_{n+2} = \left(\frac{(\wp')^2}{4(\wp+\alpha_n)(\wp+\alpha_{n+1})} - \beta_{n+1}\frac{\wp+\alpha_{n+2}}{\wp+\alpha_{n+1}} - \beta_{n}\frac{\wp+\alpha_{n-1}}{\wp+\alpha_{n}}\right)\pi_n - \beta_n\beta_{n-1}\pi_{n-2}.
\end{equation}
However, from the difference equation of the even case \eqref{differenceY2k}, we see that $\pi_{n+2}$ can be written in terms of $\pi_{n}$ and $\pi_{n-2}$ using the equation
\begin{equation}
        \pi_{n+2} = \left(\wp - B_n \right)\pi_n - \beta_n\beta_{n-1}\pi_{n-2},
\end{equation}
where $B_n$ is a constant. Moreover, the behaviour of the first term in \eqref{eq 6} near $z=0$ gives the expression
\[ B_n = \beta_{n+1}+\beta_n + \alpha_n + \alpha_{n+1}.\] 
Considering the behaviour of \eqref{eq 6} near the poles of $\wp + \alpha_n$ we further determine that
\begin{equation}\label{eq:solbetan}
    -\alpha_n^3 +\frac{g_2}{4}\alpha_n -\frac{g_3}{4} = \beta_n(\alpha_{n-1}-\alpha_n)(\alpha_{n+1}-\alpha_n),
\end{equation}
and the theorem follows.
\end{proof}

We now use the compatibility condition of the linear system \eqref{linsys:n+1}, \eqref{linsys:dz} to obtain the recurrence relation for the coefficients of the elliptic polynomials.
\begin{theorem}\label{main theorem sect 3}
The compatibility condition 
\begin{align}\label{compbegmix}
    R_n' - L_{n+1} R_n + R_n L_n = 0,
\end{align}
gives the following recurrence relations:
\begin{align}
    \beta_n&=\frac{g_3-g_2 \alpha_n+4 \alpha_n^3}{4 (\alpha_{n-1}-\alpha_n) (\alpha_n-\alpha_{n+1})},\label{rec:beta_n}\\
    \alpha_{n+1}&= \frac{(1-n) \alpha_n \left(4 \alpha_n^3-3 g_2 \alpha_n+4 g_3\right)-\alpha_{n-1} \left(4 (n-2) \alpha_n^3+n g_2 \alpha_n-(2n-1) g_3\right)}{4 n \alpha_n^3+(n-1) \alpha_{n-1} \left(g_2-12 \alpha_n^2\right)+g_2 (n-2) \alpha_n-g_2 (2 n-3)}.\label{rec:alpha_n}
\end{align}





\end{theorem}
\begin{proof}
Substituting \eqref{linsys:n+1} and \eqref{linsys:dz} into Equation \eqref{compbegmix} we find that the (1,1) entry is the only non-zero term, and is cubic in $\wp(z)$. Equating the coefficients of the powers of order 3, 2, 1, 0 to zero we find the following set of equations respectively
\begin{align}
&0=(n-2) \alpha_n-(n+1) \alpha_{n+1}+(2 n-3) \beta_n-(2 n+1) \beta_{n+1},\label{Q1}\\
  & 3 \alpha_n\alpha_{n+1}+\frac{{g_2}}{2}= ((n-2) \alpha_{n-1}+(n-1) \alpha_n+2 (2 n-3) \alpha_{n+1}) \beta_n\nonumber\\
   &\hspace{3cm}-((4 n+2) \alpha_n+n \alpha_{n+1}+(n+1) \alpha_{n+2}) \beta_{n+1}, \label{Q2}\\
    &3 g_3+ g_2 n \alpha_n-g_2 (n-1) \alpha_{n+1}= 4 \alpha_{n+1} (2 (n-2) \alpha_{n-1}+2 (n-1) \alpha_n+(2 n-3) \alpha_{n+1}) \beta_n\nonumber\\
    &\hspace{5.5cm}-4 \alpha_n ((2 n+1) \alpha_n+2 n \alpha_{n+1}+2 (n+1) \alpha_{n+2}) \beta_{n+1},\label{Q3}\\
    &-\alpha_{n+1} (g_2 \alpha_n+g_3 (n-2))+g_3 (n+1) \alpha_n =-4 (n \alpha_{n+1}+(n+1) \alpha_{n+2}) \alpha_n^2 \beta_{n+1}\nonumber\\
    &\hspace{7.5cm}+4 ((n-2) \alpha_{n-1}+(n-1) \alpha_n) \alpha_{n+1}^2 \beta_n.\label{Q4}
\end{align}
The recurrence for $\beta_n$ \eqref{rec:beta_n} follows from \eqref{eq:solbetan}. Substituting  \eqref{rec:beta_n} in \eqref{Q1}, we obtain a constraint that is cubic in $\alpha_n$, $\alpha_{n+1}$ and linear in $\alpha_{n-1}$, $\alpha_{n+2}$:
\begin{align*}
    \alpha_{n+2}=\frac{X}{Y}
\end{align*}
with
\begin{align*}
   X&= \alpha_{n+1} \left(2 g_2 (1-2 n) \alpha_n+\alpha_{n-1} \left(4 (n-2) \alpha_n^2+2 g_2 n+ g_2\right)+4 (n-1) \alpha_n^3+g_3 (2 n-3)\right)\\
   &-g_3 (2 n+1) (\alpha_{n-1}-\alpha_n)+4 n (\alpha_n-\alpha_{n-1}) \alpha_{n+1}^3-4 (2 n-1) (\alpha_{n-1}-\alpha_n) \alpha_n \alpha_{n+1}^2,\\
   Y&=\alpha_n (4 \alpha_n ((n-2) \alpha_{n-1}+(n-1) \alpha_n)+g_2 (3-2 n))+4 (n+1) (\alpha_{n-1}-\alpha_n) \alpha_{n+1}^2\\
   &+4 (2 n-1) \alpha_n (\alpha_n-\alpha_{n-1}) \alpha_{n+1}+g_3 (2 n-3).
\end{align*}
Manipulating the set of equations \eqref{Q1}-\eqref{Q4} using the above expression gives \eqref{rec:alpha_n}.
\end{proof}

Initial conditions for the recurrence \eqref{rec:alpha_n} of the $\alpha_n$ in Theorem \ref{main theorem sect 3} can be computed directly, yielding
\begin{align*}
    \alpha_3&=\frac{3g_3-4g_2\eta_1}{g_2-48\eta_1^2},\\
    \alpha_4&=\frac{5g_2^3-108g_3^2+108g_2g_3\eta_1-432g_2^2\eta_1^2+8640 g_3\eta_1^3}{18(3g_3-4g_2\eta_1)(g_2-48\eta_1^2)},
\end{align*}
where we remark that explicit formulas for $g_2,g_3$ and $\eta_1$ are given in Appendix \ref{Appendix:EllFunc}.

\section{Linear problems for the even case and the elliptic form of Painlev\'e VI}\label{Sec:Even_Lin}

In this section, we restrict to the case of the even elliptic polynomials. Firstly, we express the polynomials in terms of Hankel determinants. Secondly, the solution of the RHP yields the Lax pair of the elliptic form of Painlev\'e VI where the modular parameter $\tau$ assumes the role of the isomonodromic time. Moreover, the recurrence relation coming from the corresponding discrete linear system provides a formulation of the solution of the elliptic form of the Painlev\'e VI equation.


\subsection{Lax pair of the elliptic form of Painlev\'e VI}
\begin{theorem}\label{thm:Lax_ell_PVI}
For $k\geq 1$, the solution \eqref{rhpsol:even} solves the following pair of linear equations that correspond to the elliptic form of Painlev\'e VI 
\begin{align}
  L_{2k}(z,\tau)&:= Y'_{2k}(z,\tau) Y_{2k}(z,\tau)^{-1} = \sum_{i=1}^3 \wp'(z) \frac{  {L}_{2k}^{(i)}}{(\wp(z)-e_i)}, \label{zder:Y2k} \\
        M_{2k}(z,\tau)&:= \frac{d}{d\tau}{Y}_{2k}(z,\tau) Y_{2k}(z,\tau)^{-1} =  \sum_{i=1}^3 \frac{  {L}_{2k}^{(i)} (\dot{\wp}(z)- \dot{e}_i)}{2(\wp(z)-e_i)},   \label{dtau:M2k}
\end{align}
where the matrices $L_{2k}^{(i)}$ are given in \eqref{eq:listL2k}.
\end{theorem}
We begin with the $z$-derivative and then compute the $\tau$-derivative.
\begin{proof}

Let us the following asymptotics of the even polynomials and their Cauchy transform Lemma \ref{lemma:evenpol}\footnote{we drop the $z$, $\tau$ dependence for the remainder of the proof.}
\begin{align}
    \pi_{2n} &=\frac{1}{z^{2n}} + \frac{a_{1,2n}}{z^{2n-2}}+\frac{a_{2,2n}}{z^{2n-4}}+\mathcal{O}(z^{6-2n}), \\
    \mathcal{C}(\pi_{2n})(z) &= \frac{h_{2k}}{2\pi i}\left(z^{2n-1} + \widetilde{a}_{1,2n} z^{2n+1} +\widetilde{a}_{2,2n} z^{2n+3}+\mathcal{O}(z^{2n+5})\right).
\end{align}
    The leading behaviour of solution $Y_{2k}$ \eqref{rhpsol:even} near $z\to 0$ is then
    \begin{align}\label{asymp:Y2k4.4}
        Y_{2k}= \left(\mathbb{1} +z^2 U  + z^4 V +\mathcal{O}(z^6)  \right) \left( \begin{array}{cc}
           z^{-2k}  & 0 \\
           0  & z^{2k-3}
        \end{array} \right),
    \end{align}
    where
    \begin{align}\label{asympz:UV}
        U = \left(\begin{array}{cc}
            a_{1,2k} & \frac{h_{2k}}{2\pi i} \\
            \frac{2\pi i}{h_{2k-2}} & \widetilde{a}_{1,2k-2}
        \end{array} \right), &&  V = \left(\begin{array}{cc}
            a_{2,2k} & \frac{h_{2k}}{2\pi i}\widetilde{a}_{1,2k} \\
            \frac{2\pi i}{h_{2k-2}} a_{1,2k-2} & \widetilde{a}_{2,2k-2}
        \end{array} \right),
    \end{align}
    and due to \eqref{asymp:Cpin}, $U$, $V$ are traceless. 
    Then, in the same limit,
    \begin{align}
        Y_{2k}' Y_{2k}^{-1} \det Y_{2k} &= \frac{1}{z^4}\left( \begin{array}{cc}
           -2k  & 0 \\
           0  & 2k-3
        \end{array} \right) + \frac{1}{z^2} \left(2 U + \left[U, \left( \begin{array}{cc}
           -2k  & 0 \\
           0  & 2k-3
        \end{array} \right)  \right]  \right) \nonumber\\
        &+ 4 V - 2U^2 + \left[V, \left( \begin{array}{cc}
           -2k  & 0 \\
           0  & 2k-3
        \end{array} \right)  \right]+\mathcal{O}(z^2).
    \end{align}
    Note that the LHS of the above expression is an elliptic function
    \begin{align}
         Y_{2k}'(z+\tau) Y_{2k}(z+\tau)^{-1} = Y_{2k}'(z+1) Y_{2k}(z+1)^{-1}= Y_{2k}'(z) Y_{2k}(z)^{-1}.
    \end{align}
    Therefore using \eqref{eq:DetY2k} we have
    \begin{align}\label{zder:1}
        Y_{2k}' Y_{2k}^{-1}  = \frac{1}{\wp'(z)}\left(\widetilde{L}_{2k}^{(2)} \wp^2(z) + \widetilde{L}_{2k}^{(1)} \wp(z) +\widetilde{L}_{2k}^{(0)}\right),
    \end{align}
    where 
    \begin{align}
        \widetilde{L}_{2k}^{(2)} = -2\left( \begin{array}{cc}
           -2k  & 0 \\
           0  & 2k-3
        \end{array} \right), && \widetilde{L}_{2k}^{(1)} =-2\left(2 U + \left[U, \left( \begin{array}{cc}
           -2k  & 0 \\
           0  & 2k-3
        \end{array} \right)  \right]  \right)
    \end{align}
    \begin{align}
        \widetilde{L}_{2k}^{(0)} = -2\left(4 V - 2U^2 + \left[V, \left( \begin{array}{cc}
           -2k  & 0 \\
           0  & 2k-3
        \end{array} \right)  \right]\right),
    \end{align}
    and
    \begin{align}
        \tr  \widetilde{L}_{2k}^{(1)}  =0
        &&  \tr  \widetilde{L}_{2k}^{(0)}  = - \frac{2 h_{2k-2}}{i\pi} \left(\beta_{2k}+a_{1,2k}^2 \right).
    \end{align}
Furthermore, using the cubic equation
\begin{align}
    \left(\wp'(z) \right)^2 = 4\left(\wp -e_1 \right)\left(\wp -e_2 \right)\left(\wp -e_3 \right),
\end{align}
\eqref{zder:1} can be re-written as
\begin{align}
      L_{2k}(z,\tau):=Y_{2k}' Y_{2k}^{-1}  = \wp'(z)\left(\frac{{L}_{2k}^{(1)}}{\wp(z)-e_1}+ \frac{{L}_{2k}^{(2)}}{\wp(z)-e_2}+\frac{{L}_{2k}^{(3)}}{\wp(z)-e_3}\right),
\end{align}
with the relations 
\begin{equation}
\begin{split}\label{eq:listL2k}
    {L}_{2k}^{(1)} &= \frac{e_1 (1-e_2)  \widetilde{L}_{2k}^{(2)} +  \widetilde{L}_{2k}^{(0)}- e_1 e_3  \widetilde{L}_{2k}^{(1)}}{4(e_1-e_2)(e_1-e_3)}, \\
     {L}_{2k}^{(2)} &= \frac{e_2 (e_1 -1)  \widetilde{L}_{2k}^{(2)} + e_2 e_3  \widetilde{L}_{2k}^{(1)} -  \widetilde{L}_{2k}^{(0)}}{4(e_1 - e_2)(e_2 - e_3)}, \\
      {L}_{2k}^{(3)} &= \frac{ \widetilde{L}_{2k}^{(0)} + e_3^2  \widetilde{L}_{2k}^{(1)} + e_3 (1+ 2 e_3)  \widetilde{L}_{2k}^{(2)}}{4(e_1-e_3)(e_2-e_3)}.
\end{split}
\end{equation}

Now we derive the $\tau$-derivative.
    From \eqref{zder:Y2k} we have 
    \begin{align}
        \partial_z \log Y_{2k}(z,\tau) = \sum_{i=1}^3 {L}_{2k}^{(i)} \partial_z \log \left(\wp(z)-e_i \right),
    \end{align}
    implying that $Y_{2k}$ has the following local behaviour around $z-\frac{\omega_i}{2}\to 0$:
    \begin{align}
        Y_{2k}(z, \tau) \sim G_i \left(\wp(z)-e_i \right)^{\Lambda_i} \left(\mathbb{1}+\mathcal{O}(z^2)  \right) C_i,
    \end{align}
    where 
    \begin{align}
        \Lambda_i = G_i^{-1}   {L}_{2k}^{(i)} G_i.
    \end{align}
    The $\tau$-derivative in the vicinity of $z=\frac{w_i}{2}$ is
    \begin{align}
       & \frac{d Y_{2k}}{d\tau} = \left( \frac{d G_i}{d\tau} \left(\wp(z)-e_i \right)^{\Lambda_i} + G_{i} \Lambda_i \left( \wp(z)-e_i \right)^{\Lambda_i -\mathbb{1}} \left( \dot{\wp(z)}-\dot{e}_i \right)\right)\left(\mathbb{1}+\mathcal{O}(z^2)  \right) C_i, \nonumber\\
        \Rightarrow &\frac{d Y_{2k}}{d\tau}Y_{2k}^{-1} = \frac{d G_i}{d\tau} G_i^{-1} + \frac{  {L}_{2k}^{(i)} \left( \dot{\wp(z)}-\dot{e}_i \right)}{(\wp(z)-e_i)} + \mathcal{O}(z^4).
    \end{align}
    Moreover, near $z=0$,
    \begin{align}
        \dot{Y}_{2k} Y_{2k}^{-1} \sim z^2\left( U + \dot{U} \right)  + \mathcal{O}(z^4).
    \end{align}
Therefore, using Liouville theorem, we can uniquely determine that
    \begin{align}
         \frac{d Y_{2k}}{d\tau} Y_{2k}^{-1} =: M_{2k}(z,\tau)= \sum_{i=1}^3\frac{  {L}_{2k}^{(i)} \left( \dot{\wp}(z)-\dot{e}_i \right)}{(\wp(z)-e_i)},
    \end{align}
    and the following periodicity relations hold 
    \begin{align}
     M_{2k}(z+1,\tau)=M_{2k}(z+,\tau), &&    M_{2k}(z+\tau,\tau)=  M_{2k}(z,\tau)+  L_{2k}(z,\tau).
    \end{align}
\end{proof}
The matrices $L_{2k}$, $M_{2k}$ are in fact the Lax pair of the elliptic form of the Painlev\'e VI \cite{takasaki2001painleve}, which, with a change of variables reduces to the usual Lax pair, as will be shown in Proposition \ref{prop:PVIell_normal}. Furthermore, we can compute the eigenvalues of the residue matrices of $L_{2k}$, which turn out to have the following values owing to the specific form of the determinant of $Y_{2k}$ \eqref{eq:DetY2k}. 

For what follows it will be useful to make the linear system \eqref{zder:Y2k} traceless. To do this, we start by noting that $L_{2k}(z)$ \eqref{zder:Y2k} is an elliptic matrix function, and its trace is obtained using Jacobi's formula,
\begin{equation*}
  \tr L_{2k}(z)=\frac{d}{dz}\log \det Y_{2k}(z)=\frac{\wp''(z)}{\wp'(z)},
\end{equation*}
where used that the determinant \eqref{eq:DetY2k}
\begin{align}\label{detagain}
\det Y_{2k}(z)=-\tfrac{1}{2}\wp'(z)=f_{2k}.
\end{align}
We now use the following gauge transformation to obtain a traceless linear system:
\begin{equation}\label{tracelessgauge:Y2k}
    \mathcal{Y}_{2k}(z)=f_{2k}^{-\frac{1}{2}}Y_{2k}(z),
\end{equation}
so that
\begin{equation}\label{curlyL}
    \mathcal{Y}'_{2k}(z)=\mathcal{L}_{2k}(z)\mathcal{Y}_{2k}(z),\qquad \mathcal{L}_{2k}(z)=L_{2k}(z)-\frac{\wp''(z)}{2\wp'(z)}\mathbb{1},
\end{equation}
\begin{proposition}\label{prop:eigenvaluescurlyL2k}
    The monodromy exponents $\pm \theta_i$ around the singularities $\omega_i$, and $\pm\theta_0$ around the singularity $z=0$ of the linear system  \eqref{curlyL}, are given by
    \begin{align}
     \theta_0=\frac{1}{2}\left(4k-{3} \right), && \theta_i = \frac{1}{2}, && i=1,2,3,
    \end{align}
    and $k=1$ is the case related to self-dual Einstein metrics \cite{tod1994self, hitchin1995twistor}.
\end{proposition}
\begin{proof}
We start with the behaviour of $\mathcal{Y}_{2k}$ at $z\to 0$:
\begin{align}
    \mathcal{Y}_{2k} = (\mathbb{1}+\mathcal{O}(z)) \left(\begin{array}{cc}
        z^{-2k} & 0 \\
        0 & z^{2k-3}
    \end{array} \right) \left(\begin{array}{cc}
        z^{3/2} & 0 \\
        0 & z^{3/2}
    \end{array} \right),
\end{align}
therefore the monodromy exponents around $z=0$ are $\pm \frac{1}{2}\left(4k-{3}\right)$, which can be read off from the asymptotics \eqref{asymp:Y2k4.4}.
The coefficient matrix has simple poles at $z=w_1,w_2,w_1+w_2,0$ coming from the term $\wp'(z)$. We now determine the eigenvalues of the residue matrices at the singularities.
The asymptotic behaviour of $\mathcal{Y}_{2k}(z)$ near $z=w_i$, $i=1,2,3$ is 
\begin{equation*}
    \mathcal{Y}_{2k}(z)=c^{-\frac{1}{2}}(z-w_1)^{-\frac{1}{2}}Y_{2k}(w_1)(I+\mathcal{O}(z-w_1)),
\end{equation*}
as $f_{2k}=(z-w_i)(c+\mathcal{O}(z-w_i))$ for some $c\neq 0$. Moreover, \eqref{detagain} implies that $\det Y_{2k}(w_i)$ is zero and therefore $Y_{2k}(w_i)$ is a rank one matrix and can be expressed as
\begin{equation*}
    Y_{2k}(w_i)C=\begin{pmatrix}0 & *\\
    0 &*
    \end{pmatrix},
\end{equation*}
where the second column is nonzero with $C$ being a constant matrix. We then obtain that for $z\to w_i$,
\begin{align*}
    \mathcal{Y}_{2k}(z)C=U_0(\mathbb{1}+\mathcal{O}(z-w_i))(z-w_i)^{\frac{1}{2}\sigma_3}, && U_0\in SL_2(\mathbb{C})
\end{align*}
as $z\rightarrow w_i$, for a matrix. Therefore the monodromy exponents around $z=w_i$ are $\pm \frac{1}{2}$.
\end{proof}

We will see in what follows (from \eqref{eq:coordinates}) that the unique zero of the $(1,2)$ entry of $L_{2k}(z)$ will be $z=Q(\tau)$, where $Q(\tau)$ satisfies the elliptic form of Painlev\'e VI \cite{takasaki2001painleve} in the present case reads
\begin{align}\label{ell:PVI}
    (2\pi i)^2 \frac{d^2 Q(\tau)}{d\tau^2}=\sum_{i=0}^3 \alpha_i \wp'(Q(\tau)+w_i),
\end{align}
where $w_0=0$, and
\begin{align}
  \alpha_0=  \frac{(\theta_0-1)^2 }{2}= \frac{(4k-5)^2}{8}, && \alpha_1=-\frac{\theta_1^2}{2}=-\frac{1}{8},&& \alpha_2=\frac{\theta_3^2}{2}=\frac{1}{8},&& \alpha_3=\frac{\left(1-\theta_2^2 \right)}{2}=\frac{3}{8}.
\end{align}
For $k=1$ this is the Hitchin case \cite{hitchin1995twistor, manin1998universal}, and the corresponding solution $Q=Q(\tau)$, is given by
\begin{equation}\label{eq:hitchin_case}
     \wp(Q(\tau))=e_1+e_2+2\frac{16\eta_1^3+g_2\eta_1-g_3}{48\eta_1^2-g_2}.
\end{equation}
This formula follows from equation \eqref{eq:uk_gamma_delta}.

We finish this subsection with a recursive formula for the $Y_{2k}$.
\begin{proposition}\label{prop:rhprec}
The matrix function $Y=Y_{2k}$ satisfies the following discrete evolution with respect to $k$,
    \begin{align}\label{differenceY2k}
        Y_{2k+2} = R_{2k} Y_{2k}, && R_{2k}:=\left(\begin{array}{cc}
         \wp(z)+a_{1,2k+2}- a_{1,2k}   & -\frac{h_{2k}}{2\pi i} \\
          \frac{2\pi i}{h_{2k}}   & 0
        \end{array} \right),
    \end{align}
    for $k\geq 1$.
\end{proposition}
\begin{proof}
    The proof is similar to that of Proposition \ref{prop:difference_ls_Yn} and follows from a direct substitution of \eqref{rhpsol:even} with the expressions of Lemmas \ref{lemma:evenpol} and \ref{lemma:detY2k}.
\end{proof}
From the above proposition, we obtain the following three-term recurrence relation for the orthogonal polynomials,
\begin{align}
   \wp(z) \pi_{2k} =  \pi_{2k+2} + \widetilde{\alpha}_{2k} \pi_{2k} + \beta_{2k} \pi_{2k-2},
\end{align}
which coincides with \eqref{eq:rec_pin2pin}, as can be seen by using the relations \eqref{def:alphabetan}, \eqref{asymp:Cpinmixed}, and \eqref{asymp:Cpin}.

\subsection{Relation to Painlev\'e VI}
Under a simultaneous transformation of the dependent and independent variables,
\begin{align}\label{change of variable:tut}
    t = \frac{e_3-e_1}{e_2-e_1}, && u(t) = \frac{\wp(Q(\tau))-e_1}{e_2-e_1},
\end{align}
the elliptic form of Painlev\'e VI \eqref{ell:PVI} reduces to the usual form \cite{takasaki2001painleve}
\begin{align*}
			\ddot{u}=&\left(\frac{1}{u}+\frac{1}{u-1}+\frac{1}{u-t}\right)\frac{\dot{u}^2}{2}-\left(\frac{1}{t}+\frac{1}{t-1}+\frac{1}{u-t}\right)\dot{u}\\
			&+\frac{u(u-1)(u-t)}{2t^2(t-1)^2}\left((2k-\tfrac{1}{2})^2-\frac{ t}{4u^2}+\frac{(t-1)}{4(u-1)^2}+\frac{3t(t-1)}{4(u-t)^2}\right),
\end{align*}
where $u\equiv u_{k}(t)$, and $\dot{}$ denotes derivative with respect to $t$.

In this section, we derive a corresponding transformation
for the Lax pair.

\begin{proposition}\label{prop:PVIell_normal}
The change of variables
    \begin{align}\label{ell_sph:cov}
        x= \frac{\wp(z)-e_1}{e_2-e_1}, && t = \frac{e_3-e_1}{e_2-e_1}, && \widetilde{Y}_{2k}(x,t)=\mathcal{Y}_{2k}(z,\tau)\left(  e_1-e_2\right)^{-(k-\frac{3}{4})\;\sigma_3},
    \end{align}
transforms the linear system in Theorem \ref{thm:Lax_ell_PVI} into the following  $4$-point Fuchsian system, with the singularities at $w_1, w_2, w_3, 0$ mapped to $0,1,t, \infty$ respectively, and corresponding deformation equation,
\begin{equation}\label{eq:pvi_lax}
  \frac{d\widetilde{Y}_{2k}}{dx}= A_{2k} \widetilde{Y}_{2k},\quad \frac{d\widetilde{Y}_{2k}}{dt}= B_{2k} \widetilde{Y}_{2k},
\end{equation}    
with
    \begin{align*}
        A_{2k}(x,t)= \frac{A_{2k}^{(1)}}{x} + \frac{A_{2k}^{(2)}}{x-1} + \frac{A_{2k}^{(3)}}{x-t}, \qquad B_{2k}(x,t) = - \frac{A_{2k}^{(3)}}{x-t},
    \end{align*}
 and coefficient matrices above are related to $\mathcal{L}_{2k}$ as
 \begin{align}\label{gauge:AL}
    A^{(i)}_{2k}= \left(  e_2-e_1\right)^{(4k-3)\sigma_3/4} \mathcal{L}_{2k}^{(i)} \left(  e_2-e_1\right)^{-(4k-3)\sigma_3/4}, \qquad i=1,2,3.
\end{align}   
In particular, the monodromy of $\widetilde{Y}_{2k}(x,t)$, with respect to $x$, is constant in $t$.
\end{proposition}
\begin{proof}
   The linear system \eqref{curlyL}:
    \begin{align*}
    \frac{\partial}{\partial z}\mathcal{Y}_{2k}(z,\tau) = \wp'(z)\left(\frac{\mathcal{L}_{2k}^{(1)}(\tau)}{\wp(z)-e_1}+\frac{\mathcal{L}_{2k}^{(2)}(\tau)}{\wp(z)-e_2} +\frac{\mathcal{L}_{2k}^{(3)}(\tau)}{\wp(z)-e_3}  \right) \mathcal{Y}_{2k}(z,\tau)
    \end{align*}
   under the change of  variables \eqref{ell_sph:cov} reads
    \begin{align}\label{eq:lin_transf_prop5}
    \frac{\partial}{\partial x}\mathcal{Y}_{2k}(x,t) = \left(\frac{\mathcal{L}_{2k}^{(1)}(t)}{x}+\frac{\mathcal{L}_{2k}^{(2)}(t)}{x-1} +\frac{\mathcal{L}_{2k}^{(3)}(t)}{x-t}  \right) \mathcal{Y}_{2k}(x,t).
    \end{align}
Let us now understand the change of variables for the equation \eqref{dtau:M2k} $\tau$ after the gauge transformation \eqref{tracelessgauge:Y2k}:
    \begin{align}\label{eq:ODE_Y_tau}
        \frac{d}{d\tau} \mathcal{Y}_{2k}(z,\tau)  &=  \left(\sum_{i=1}^3 \frac{(\dot{\wp}(z)-\dot{e}_i) \mathcal{L}_{2k}^{(i)}}{(\wp(z) -e_i)}\right) \mathcal{Y}_{2k}(z,\tau).
    \end{align}
We begin by simplifying the sum on the right-hand side. To this end,
    we note that the definition of $x$ in equation \eqref{ell_sph:cov}, and the fact that $\frac{dx }{d\tau} =  0$, imply
\begin{align*}
   \dot{\wp}(z) &= \frac{\left(\wp(z)-e_1 \right) \dot{e}_2 - \left(\wp(z)-e_2 \right) \dot{e}_1}{e_2 - e_1}.
\end{align*}
From the above equation, and the identity $e_1+e_2+e_3=0$, we obtain the following relations
\begin{align*}
\frac{\dot{\wp}(z)- \dot{e}_1}{\wp(z) - e_1} &= \frac{\dot{e}_2 - \dot{e}_1 }{e_2 - e_1} , \\
\frac{\dot{\wp}(z)- \dot{e}_2}{\wp(z) - e_2}& = \frac{\dot{e}_2 - \dot{e}_1 }{e_2 - e_1} ,\\
\frac{\dot{\wp}(z)- \dot{e}_3}{\wp(z) - e_3} &=
\frac{\dot{e}_2 - \dot{e}_1 }{e_2 - e_1}
-3\frac{e_1 \dot{e}_2 -\dot{e}_1 e_2}{(e_2 - e_1)(\wp(z) - e_3)}.
\end{align*}
Using the above identities, we can express the sum on the right-hand side of equation \eqref{eq:ODE_Y_tau} rationally in $x$, giving
\begin{align}\label{eq:onlyxtransf}
    \frac{d}{d\tau} \mathcal{Y}_{2k}(x,\tau) \mathcal{Y}_{2k}(x,\tau)^{-1}=   \left(\mathcal{L}_{2k}^{(1)}+\mathcal{L}_{2k}^{(2)}+\mathcal{L}_{2k}^{(3)} \right)\frac{(\dot{e}_2 - \dot{e}_1)}{e_2 - e_1} - \frac{3 \mathcal{L}_{2k}^{(3)} \left(e_1 \dot{e}_2 - e_2 \dot{e}_1 \right)}{\left(x-\left( \frac{e_3-e_1}{e_2-e_1}\right)\right) (e_2-e_1)^2}, 
\end{align}
where 
\begin{align}
\mathcal{L}_{2k}^{(1)}+\mathcal{L}_{2k}^{(2)}+\mathcal{L}_{2k}^{(3)}=\frac{1}{4} (4k-3)\sigma_3,
\end{align}
as can be seen from \eqref{eq:listL2k}, along with the gauge transformation \eqref{curlyL}.
We now observe that, the gauge transformation
\begin{align}\label{gauge:tildeY}
    \widetilde{{Y}}:= \mathcal{Y}_{2k} \left(  e_1-e_2\right)^{-(4k-3)\;\sigma_3/4},
\end{align}
removes the constant term in \eqref{eq:onlyxtransf}, giving
\begin{align}
     \frac{d}{d\tau} \widetilde{Y}(x,\tau) \widetilde{Y}(x,\tau)^{-1}=   - \frac{3 \left(e_1 \dot{e}_2 - e_2 \dot{e}_1 \right)}{(e_2-e_1)^2}\frac{A^{(3)}_{2k}}{\left(x-\left( \frac{e_3-e_1}{e_2-e_1}\right)\right)}. 
\end{align}
Using the identity below coming from the change of variable $\tau\to t$
\begin{align}
      \frac{dt}{d\tau} = 3 \frac{(\dot{e}_2 e_1 - \dot{e}_1 e_2 )}{(e_2-e_1)^2}
\end{align}
implies that
\begin{align}
     \frac{d}{dt} \widetilde{Y}(x,t) \widetilde{Y}(x,t)^{-1}=   - \frac{A^{(3)}_{2k}(t)}{\left(x-t\right)}.
\end{align}
Furthermore, \eqref{eq:lin_transf_prop5} with the gauge transformation \eqref{gauge:tildeY} and \eqref{gauge:AL} is
\begin{align}\label{eq:AinY}
   \frac{d}{dx} \widetilde{Y}(x,t) \widetilde{Y}(x,t)^{-1} = \frac{A^{(1)}_{2k}}{x} + \frac{A^{(2)}_{2k}}{x-1}+ \frac{A^{(3)}_{2k}}{x-t}. 
\end{align}
which finishes the proof of the proposition.
\end{proof}

We note that the following rational matrix is traceless,
\begin{equation*}
    A^{(1)}_{2k}+A^{(2)}_{2k}+A^{(3)}_{2k}=-\frac{\theta_0}{2}\sigma_3,\quad \theta_0:=\tfrac{3}{2}-2k,
\end{equation*}
and
\begin{equation*}
    |A^{(i)}_{2k}|=-\frac{\theta_i^2}{4},\quad \theta_i\equiv \theta=\frac{1}{2}\qquad (i=1,2,3),
\end{equation*}
as can be obtained from Proposition \ref{prop:eigenvaluescurlyL2k} and \eqref{gauge:AL}.
Introducing standard coordinates $(u_k,v_k,g_k)$, see e.g. \cite{jimbomiwaII}, through
\begin{align}
    &(A_{2k})_{12}(x,t)=-\theta_0 g_k\frac{x-u_k}{2 x(x-t)(x-1)},\label{eq:coordinates}\\
    &(A_{2k})_{11}(x,t)=v_k-\frac{\theta}{2}\left(\frac{1}{u_k}+\frac{1}{u_k-t}+\frac{1}{u_k-1}\right), \nonumber
\end{align}
we obtain that $u\equiv u_k$ and $v\equiv v_k$ satisfy the equations below
\begin{align}
    t(t-1)\dot{u}&=-\tfrac{1}{2}t+u\left(t-\tfrac{1}{2}u+2p(u-t)(u-1)\right),\label{eq:yt}\\
    t(t-1)\dot{v}&=\tfrac{3}{16}+\frac{\theta_0}{4}(\theta_0-2)-t v(1+v)+ v u+2(1+t)v^2u-3v^2u^2,
\end{align}
and $g(t)=g_k(t)$ solves
\begin{equation}\label{eq:gauge_differential}
    \frac{g_k'(t)}{g_k(t)}=(1-4k)\frac{u(t)-t}{2 t(t-1)}.
\end{equation}
Note that the zero of $\left(\mathcal{L}_{2k}\right)_{12}$ being $z=Q(\tau)$ follows from \eqref{eq:coordinates} under the change of variable \eqref{change of variable:tut}.

\section{The Painlev\'e VI tau-function and Hankel determinants}\label{Sec:HD_PVI}
In this section, we find explicit formulas for the solution $u_k$ of Painlev\'e VI, introduced in \eqref{eq:coordinates}, in terms of the Hankel determinants of moments defined in equation \eqref{eq:delta2k}. To this end, we compute the first couple of matrix coefficients in the the expansion around $z=\infty$ of the solution to RHP \ref{rhp:even}.

\begin{lemma}
The asymptotic expansion of the solution of RHP \ref{rhp:even} around $z=0$ can be written as
\begin{equation}\label{eq:asymp_Y2k}
    Y_{2k}(z) = \left(\mathbb{1} + \wp(z)^{-1} U + \wp(z)^{-2} V + \mathcal{O}(\wp(z)^{-3}) \right)\left(\begin{array}{cc}
      \wp(z)^k   & 0 \\
        0 & \wp(z)^{-k}
    \end{array} \right) \left(\begin{array}{cc}
      1   & 0 \\
        0 & -\frac{1}{2}\wp'(z)
    \end{array} \right),
\end{equation}
where the matrices $U$ and $V$ are given explicitly by
\begin{equation*}
    U=\begin{pmatrix}
    -\frac{\Gamma_{2k}}{\Delta_{2k}} & \frac{\Delta_{2k+2}}{2\pi i \Delta_{2k}}&\\
    \frac{2\pi i\Delta_{2k-2}}{\Delta_{2k}} & +\frac{\Gamma_{2k}}{\Delta_{2k}}
    \end{pmatrix},\qquad
    V=\begin{pmatrix}\vspace{1mm}
    \frac{\Lambda_{2k}}{\Delta_{2k}} & \frac{\Gamma_{2k+2}}{2\pi i \Delta_{2k}}&\\
    -\frac{2\pi i\Gamma_{2k-2}}{\Delta_{2k}} & {\sf v}_{22}
    \end{pmatrix},
\end{equation*}
with
\begin{equation*}
    {\sf v}_{22}=\frac{\Gamma_{2k}^2}{\Delta_{2k}^2}+\frac{\Delta_{2k-2}\Delta_{2k+2}}{\Delta_{2k}}-\frac{\Lambda_{2k}}{\Delta_{2k}}.
\end{equation*}
\end{lemma}
\begin{proof} Recall the explicit, and unique, solution $Y_{2k}(z)$ of RHP \ref{rhp:even}, defined in equation \eqref{rhpsol:even}. Note that $\widehat{Y}(z)=Y_{2k}(-z)\sigma_3$, also satisfies all the conditions in the RHP, and thus
\begin{equation*}
    Y_{2k}(z)=Y_{2k}(-z)\sigma_3.
\end{equation*}
It follows from this symmetry that $Y_{2k}(z)$ admits an expansion in powers of $\wp(z)$ as given in the lemma.

We proceed to compute the coefficient matrices $U$ and $V$. The expressions for ${\sf u}_{11}$, ${\sf v}_{11}$, ${\sf u}_{21}$ and ${\sf v}_{21}$ follow directly from the expansions of the corresponding orthogonal polynomials in equation \eqref{eq:OPexpansion}.
Next, by equation \eqref{asymp:Cpinmixed}, we have
\begin{equation*}
   \mathcal{C}(\pi_{2k})(z)=\frac{h_{2k}}{2\pi i} z^{2k-1}(1+\mathcal{O}(z^2)),
\end{equation*}
from which the expression for $u_{12}$ follows. Finally, note that $|Y_{2k}(z)|=-\tfrac{1}{2}\wp'(z)$ implies
\begin{equation*}
|\mathbb{1}+\wp(z)^{-1}U+\wp(z)^{-2}V|=1+\mathcal{O}(\wp(z)^{-3})
\end{equation*}
as $z\rightarrow 0$, which is equivalent to
\begin{equation*}
    \operatorname{Tr} U=0,\quad \operatorname{Tr}  V+|U|=0,
\end{equation*}
as can also be seen from \eqref{asymp:Y2k4.4}. \eqref{asympz:UV}.
The expressions for ${\sf u}_{22}$ and ${\sf v}_{22}$ are obtained from the above two equations.
\end{proof}

It follows from the asymptotic expansion \eqref{eq:asymp_Y2k} for $Y_{2k}(z)$, that $\widetilde{Y}_{2k}(x)$ has an expansion around $x=\infty$ of the form
\begin{align*}
\widetilde{Y}_{2k}(x)&=\Psi_{2k}(x)G(x),\\
\Psi_{2k}(x)&=I+x^{-1}\widetilde{U}+x^{-2}\widetilde{V}+\mathcal{O}(x^{-3}),\\
G(x):&=x^{(k-\frac{3}{4}) \sigma_3}\bigg(1+\frac{e_2}{(e_1-e_2)x}\bigg)^{k\;\sigma_3}((1-t/x)(1-1/x))^{-\frac{1}{4}\sigma_3},
\end{align*}
where
\begin{equation*}
    \widetilde{U}=\frac{U}{e_1-e_2},\qquad \widetilde{V}=\frac{V-e_2\;U}{(e_1-e_2)^2}.
\end{equation*}
Now, by equation \eqref{eq:AinY}, we can express the coefficient matrix $A_{2k}(x)$, in terms of $\Psi_{2k}(x)$, as follows,
\begin{equation}\label{eq:A2kequation}
\widetilde{A}_{2k}(x)=\Psi_{2k}'(x)\Psi_{2k}(x)^{-1}+
\Psi_{2k}(x)G'(x)G(x)^{-1}\Psi_{2k}(x)^{-1}.
\end{equation}
This expression allows us to compute the coordinates $(u_k,v_k,g_k)$ introduced in equation \eqref{eq:coordinates}. Indeed, by expanding the right-hand side of equation \eqref{eq:A2kequation} around $x=\infty$, we get the following expansion for its $(1,2)$-entry,
\begin{align*}
 (\widetilde{A}_{2k})_{12}(x)&=-\tfrac{1}{2}(4k-1)\widetilde{\sf u}_{12}x^{-2}+\tfrac{1}{2}s\;x^{-3}+\mathcal{O}(x^{-4}),\\
 s:&=-(4k+1)t-(4k-3)\widetilde{\sf u}_{11}\widetilde{\sf u}_{12}+\left(1+t+4k\frac{e_2}{e_1-e_2}\right)\widetilde{\sf u}_{12}.
\end{align*}
Comparing this asymptotic expansion with equation \eqref{eq:coordinates}, and recalling that
\begin{equation*}
    e_1-e_2=4\mathcal{K}(t)^2,
\end{equation*}
where $\mathcal{K}(t)$ is the complete elliptic integral of the first kind (see Appendix \ref{Appendix:EllFunc}),
we obtain
\begin{align}\label{eq:uk_gamma_delta}
    u_k(t)&=\frac{1}{4\mathcal{K}(t)^2}\left(\frac{2k-3}{2k-1}\frac{\Gamma_{2k}}{\Delta_{2k}}-\frac{2k+1}{2k-1}\frac{\Gamma_{2k+1}}{\Delta_{2k+1}}\right)+\frac{1+t}{3},\\
    g_k(t)&=-\frac{(2k-1)}{4\pi i(2k-3)}\mathcal{K}(t)^{1-2k}h_{2k}(t).\label{eq:gauge_factor}
\end{align}

We are now in a position to prove the following theorem.
\begin{theorem}\label{thm:pvisolution} For $k\geq 0$,
    \begin{equation}\label{eq:y_formula}
    u_k(t)=\frac{2t(t-1)}{2k-1}\left(\frac{\dot{\Delta}_{2k}(t)}{{\Delta}_{2k}(t)}-\frac{\dot{\Delta}_{2k+2}(t)}{{\Delta}_{2k+2}(t)}\right)+1-\frac{\mathcal{E}(t)}{\mathcal{K}(t)},
\end{equation}
solves $P_\text{VI}$, where $\mathcal{K}(t)$ and $\mathcal{E}(t)$ denote the complete elliptic integrals of the first and second kind respectively (see Appendix \ref{Appendix:EllFunc}), with parameter values
\begin{equation*}
\theta_1=\theta_2=\theta_3=\tfrac{1}{2},\quad 
    \theta_0=\tfrac{1}{2}(3-2k).
\end{equation*}
\end{theorem}
\begin{proof}
Recalling that $\Delta_0(t)=\Delta_2(t)=1$ and
\begin{equation*}
    \Delta_4(t)=\tfrac{16}{3}\mathcal{K}(t)^2((t-1)\mathcal{K}(t)^2-2(t-2)\mathcal{K}(t)\mathcal{E}(t)-3\mathcal{E}(t)^2),
\end{equation*}
it can be check by direct calculation that
\begin{subequations}\label{eq:y0y1}
\begin{align}
    u_0(t)&=1-\frac{\mathcal{E}(t)}{\mathcal{K}(t)},\\
    u_1(t)&=1-\frac{\mathcal{E}(t)}{\mathcal{K}(t)}-\frac{2t(t-1)}{3}\frac{\dot{\Delta}_{4}(t)}{\Delta_{4}(t)},
\end{align}
\end{subequations}
solve $P_\text{VI}$ for the parameter values indicated in the theorem.

Now, assume $k\geq 2$. By combining the differential equation for the gauge factor $g_k(t)$, equation \eqref{eq:gauge_differential}, with the explicit expression for $g_k(t)$ in terms of $h_{2k}(t)$, equation \eqref{eq:gauge_factor}, we obtain
\begin{equation*}
    \frac{\dot{h}_{2k}(t)}{h_{2k}(t)}=-\frac{(2k-1)}{2t(t-1)}\left(\frac{\mathcal{E}(t)}{\mathcal{K}(t)}+u_k(t)-1\right).
\end{equation*}
Solving this equation for $u_k(t)$ and using equation \eqref{eq:normdelta}, we obtain the expression for $u_k(t)$ given in the theorem. Since we already know that $u_k(t)$ solves Painlev\'e VI, the theorem follows.
\end{proof}

\subsection{The Painlev\'e VI $\tau$-function}
The Painlev\'e VI tau function $\mathcal{T}_k(t)$, corresponding to the linear system \eqref{eq:pvi_lax},
can be defined by
\begin{equation}\label{eq:tau_definition}
    \zeta_k(t)=t(t-1)\frac{d}{d t}\log{\mathcal{T}_k(t)},
\end{equation}
up to a multiplicative constant, where
\begin{align}
    \zeta_k=&(t-1)\operatorname{Tr}{{A}_2{A}_3}+t \operatorname{Tr}{{A}_1{A}_3}\nonumber\\
    =&u_k(u_k-t)(u_k-1)v_k^2-\tfrac{1}{2}(t-2(1+t)u_k+3u_k^2)v_k\label{eq:zeta}\\
    &-\tfrac{1}{2}k(2k-3)u_k+\tfrac{1}{4}(4k^2-6k+1)t-\tfrac{1}{8}.\nonumber
\end{align}
It is an analytic function on the universal covering space of the punctured sphere $\mathbb{CP}^1\setminus\{0,1,\infty\}$, and satisfies the ODE
 \begin{equation*}
     (t(t-1)\ddot{\zeta}(t))^2=-2\begin{vmatrix}
     \tfrac{1}{8} & t\dot{\zeta}-\zeta & \tfrac{3}{8} -\frac{1}{4}\theta_0^2+\dot{\zeta}\\
     t\dot{\zeta}-\zeta &  \tfrac{1}{8}   & (t-1)\dot{\zeta}-\zeta\\
     \tfrac{3}{8} -\frac{1}{4}\theta_0^2+\dot{\zeta} &
    (t-1)\dot{\zeta}-\zeta & \tfrac{1}{8}
     \end{vmatrix}.
 \end{equation*}

In the following theorem, we give an explicit expression for $\mathcal{T}_k(t)$.
\begin{theorem}\label{thm:tau} For $k\geq 0$,
    \begin{equation}\label{eq:tau_delta}
\mathcal{T}_k(t)=t^{\frac{1}{8}}(1-t)^{\frac{1}{8}}(2\mathcal{K}(t))^{-n(2k-3)}\Delta_{2k}(t).
\end{equation}
\end{theorem}
\begin{proof}
To prove the theorem, it is enough to derive the following expression for $\zeta_k(t)$,
\begin{equation}\label{eq:zetan}
    \zeta_k(t)=t(t-1)\frac{\dot{\Delta}_{2k}(t)}{\Delta_{2k}(t)}+\frac{1}{2}(2k)(2k-3)\left(\frac{\mathcal{E}(t)}{\mathcal{K}(t)}+t-1\right)+\frac{1}{8}(2t-1).
\end{equation}
for $k\geq 0$. We will prove this expression by induction.

We first deal with the cases $k=0$ and $k=1$. We recall the explicit expressions for $u_0$ and $u_1$ in equations \eqref{eq:y0y1}. Using equation \eqref{eq:yt}, we obtain the following corresponding expressions for $v_0$ and $v_1$,
\begin{align*}
    v_0(t)&=0,\\
    v_1(t)&=-\frac{3\Delta_4(t)}{32\mathcal{K}(t)\mathcal{E}(t)(\mathcal{K}(t)-\mathcal{E}(t))((1-t)\mathcal{K}(t)-\mathcal{E}(t)},
\end{align*}
and consequently, using equation \eqref{eq:zeta}, we obtain the following expressions for $\zeta_0(t)$ and $\zeta_1(t)$,
\begin{align*}
    \zeta_0(t)&=\tfrac{1}{8}(2t-1),\\
    \zeta_1(t)&=\tfrac{1}{8}(2t-1)-\tfrac{1}{2}\left(\frac{\mathcal{E}(t)}{\mathcal{K}(t)}+t-1\right).
\end{align*}
This shows that equation \eqref{eq:zetan} holds for $k=0,1$.

Next, we derive a recursive equation for $\zeta_k(t)$. To this end, we note that the recurrence in Proposition \ref{prop:rhprec}, translates to the following recurrence for $\widetilde{Y}_{2k}(x,t)$,
\begin{equation*}
    \widetilde{Y}_{2k+2}(x,t)=R_{2k}(x,t)\widetilde{Y}_{2k}(x,t),
\end{equation*}
where
\begin{equation*}
    R_{2k}(x,t)=\begin{pmatrix}x+r_{11} & r_{12}\\
r_{21} & 0
    \end{pmatrix},
\end{equation*}
with
\begin{equation*}
    r_{12}=\frac{\theta_0 k}{2(\theta_0-1)},\quad r_{21}=-\frac{2(\theta_0-1)}{\theta_0 k},
\end{equation*}
and
\begin{align*}
    \frac{1}{4}\theta_0(\theta_0-2)r_{11}=&\tfrac{1}{2}u_k(u_k-t)(u_k-1)v_k^2-\tfrac{1}{4}\left(3u_k^2-2(t+1)u_k+t\right)v_k\\
    &+\tfrac{3}{32}(3u_k-t-1)+\tfrac{1}{8}\theta_0(\theta_0-2)(u_k-t-1).
\end{align*}
Here $\theta_0=\tfrac{3}{2}-k$, as before.

From this recursive formula, we obtain the following recurrence for the coefficient matrix of the linear system in equation \eqref{eq:pvi_lax},
\begin{equation*}
\widetilde{A}_{2k+2}(x,t)=R_{2k}(x,t)\widetilde{A}_{2k}(x,t) R_{2k}(x,t)^{-1}+\left(\frac{\partial}{\partial x}R_{2k}(x,t)\right) R_{2k}(x,t)^{-1}.
\end{equation*}

Direct substitution now yields a very compact recursive formula for $\zeta_k(t)$,
\begin{equation}\label{eq:zeta_rec}
   \zeta_{k+1}(t)=\zeta_{k}(t)+(\theta_0-1)(u_k(t)-t).
\end{equation}
By combining this recursive formula with the equation for $u_k(t)$ in Theorem \ref{thm:pvisolution}, equation \eqref{eq:zetan} follows by induction. This completes the proof of the theorem.
\end{proof}

\begin{corollary}
From equations \eqref{eq:zeta_rec}, \eqref{eq:y_formula} and \eqref{eq:tau_definition}, we obtain the following recursion for the Painlev\'e VI tau function $\mathcal{T}_{k}(t)$,
\begin{align*}
s_k\mathcal{T}_{k-1}(t)\mathcal{T}_{k+1}(t)=&4(4k-3)^2t^2(t-1)^2\mathcal{T}_{k}(t)\ddot{\mathcal{T}}_{k}(t)-4(4k-1)(4k-5)t^2(t-1)^2 \dot{\mathcal{T}}_{k}(t)^2\\
&+2((4k-3)^2+1)t(t-1)(2t-1)\mathcal{T}_{k}(t)\dot{\mathcal{T}}_{k}(t)\\
&+\left[2(k-1)(2k-1)(4k^2-6k+1+t-t^2)-\tfrac{1}{4}\right]\mathcal{T}_{k}(t)^2,
\end{align*}

for $k\geq 1$. Here, the $s_k$ are some nonzero constants which are not rigidly defined in general. However, using the exact formula \eqref{eq:tau_delta}, they become numerical constants, and the first few are given by
\begin{equation*}
    s_1=-3,\quad s_2=525,\quad s_3=6237,\quad s_4=27885,\quad s_5=82365.
\end{equation*}
\end{corollary}



\appendix
\addtocontents{toc}{\protect\setcounter{tocdepth}{0}}

\section{Elliptic functions and their periodicity properties}\label{Appendix:EllFunc}

The Weierstrass cubic reads as
\begin{align}
    \left(\wp'(z)\right)^2 = 4 \wp^3(z) - g_2 \wp(z)- g_3 = 4 \left(\wp(z)- e_1 \right)\left(\wp(z)- e_2 \right)\left(\wp(z)- e_3 \right),
\end{align}
where 
\begin{align}
    e_1 = \wp\left( \frac{1}{2}\right), &&   e_2 = \wp\left( \frac{\tau}{2}\right), &&   e_3 = \wp\left( \frac{1+\tau}{2}\right).
\end{align}
The Weierstrass $\wp$-function is doubly periodic
\begin{align}
    \wp(z+1) =\wp(z), && \wp(z+\tau) =\wp(z),
\end{align}
and has a double pole at zero
\begin{align}
    \lim_{z\to 0} z^2 \wp(z) =1.
\end{align}
The Weierstrass $\zeta$-function is defined to be the anti-derivative of $\wp(z)$ uniquely characterised by
\begin{align}
    \zeta'(z) = -\wp(z),\quad \zeta(z)=\frac{1}{z}+\mathcal{O}(z)\quad (z\rightarrow 0),
\end{align}
and has the following periodic properties
\begin{align}\label{periodicity:zeta}
    \zeta(z+ 1) = \zeta(z)+ \eta_1(\tau), && \zeta(z+ \tau) = \zeta(z)+ \eta_2(\tau),
\end{align}
which in turn define the Weierstrass $\eta$-functions.

The elliptic nome is defined by
\begin{equation*}
    q=\exp{\frac{i\pi \omega_3}{\omega_1}}=\exp{i\pi(1+\tau)},
\end{equation*}
and we define
\begin{equation*}
    t=\frac{e_3-e_2}{e_1-e_2}=\lambda(\tau),
\end{equation*}
where $\lambda(\cdot)$ is the modular lambda function.

We have the following explicit expressions for $(e_1,e_2,e_3)$ in terms of $q$ and $t$,
\begin{align*}
(e_1,e_2,e_3)=&\frac{\pi^2}{3}\left(\theta_3(0,q)^4+\theta_4(0,q)^4,\theta_2(0,q)^4-\theta_4(0,q)^4),-\theta_2(0,q)^4-\theta_3(0,q)^4\right)\\
    =&\frac{4}{3}\mathcal{K}(t)^2(2-t,-1-t,2t+1),
\end{align*}
where $\theta_j(z,q)$ denotes the $j$th Jacobi elliptic function for $1\leq j\leq 4$. In particular
\begin{equation*}
    \theta_2(0,q)=2q^{\frac{1}{4}}\sum_{n=0}^\infty q^{n(n+1)},\quad
    \theta_3(0,q)=1+2\sum_{n=0}^\infty q^{n^2},\quad
    \theta_4(0,q)=\theta_3(0,-q).
\end{equation*}

This yields the following formulas for the invariants $\{g_2,g_3\}$ in terms of $q$ and $t$,
\begin{align*}
    g_2&=\frac{64}{3}(t^2-t+1)\mathcal{K}(t)^4,\\
    &=\frac{4\pi^4}{3}(\theta_2(0,q)^8-\theta_2(0,q)^4\theta_3(0,q)^4+\theta_3(0,q)^8),\\
    g_3&=\frac{256}{27}(2t-1)(t-2)(t+1)\mathcal{K}(t)^6,\\
    &=\frac{8\pi^6}{27}\left(\theta_2(0,q)^{12}+\theta_3(0,q)^{12}-\frac{3}{2}\theta_2(0,q)^4\theta_3(0,q)^4(\theta_2(0,q)^4+\theta_3(0,q)^4\right).
\end{align*}

Finally, we note the following useful formula for $\eta_1$,
\begin{align*}
    \eta_1=&-\frac{\pi^2}{6}\frac{\theta_1^{(3)}(0,q)}{\theta_1^{(1)}(0,q)}\\
    =&\frac{2}{3}\mathcal{K}(t)((t-2)\mathcal{K}(t)+3\mathcal{E}(t)),
\end{align*}
where $\theta_1^{(j)}(z,q)$ denotes the $j$th derivative of $\theta_1(z,q)$ with respect to $z$.



\section{List of polynomials}

\begin{align*}
    \pi_0 &= 1,\\
    \pi_2 &= \wp(z)+a_{1,1},\\
    \pi_3 &= - \tfrac{1}{2}\wp'(z),\\
    \pi_4 &= \wp^2(z) +a_{1,2} \wp +a_{2,2},\\
    \pi_5 & = - \tfrac{1}{2}\wp'(z) \wp + b_{2,2} \wp(z).
\end{align*}
We compute the first few coefficients
\begin{align}
    a_{1,1} = 2 \eta_1(\tau), \,\, a_{1,2}  = 2\left(\frac{4\eta_1(\tau) g_2-3g_3}{5 g_2 - 240 \eta_1^2(\tau)} \right), \,\, a_{2,2} = 4\eta_1(\tau) \left(\frac{4\eta_1(\tau) g_2-3g_3}{5 g_2 - 240 \eta_1^2(\tau)} \right) -\frac{g_2}{12},
\end{align}
and 
\begin{align}
    h_1(\tau)= 1,&& h_2(\tau)= -4\eta_1^2(\tau)+\frac{g_2}{12}, && h_3(\tau) = \frac{1}{5} \left(g_2 \eta_1^2(\tau)- 3 g_3 \right).
\end{align}

\section{Structure of the moment matrix}\label{appendix:mom}
Let us elaboratte on the structure of the moment matrix $D_n$ \eqref{eq:Dn}. There are three points to note about the moments \eqref{def:moment}:
\begin{enumerate}
    \item all the mixed moments vanish, {\it i.e} for all $i, j$
    \begin{align}\label{propmom:1}
        \mu_{2i, 2j+1} =0, 
    \end{align}
    \item the following symmetry property holds 
    \begin{align}\label{propmom:2}
    \mu_{i,j}= \mu_{j,i},
    \end{align}
    \item and, generally \begin{align}\mu_{2i, 2j}\neq \mu_{2i-1,2j+1}.\end{align} 
\end{enumerate}
With the above properties, following figure illustrates the moment matrix for $n=9$.
\begin{figure}[H]
    \centering
    \includegraphics[width=11cm]{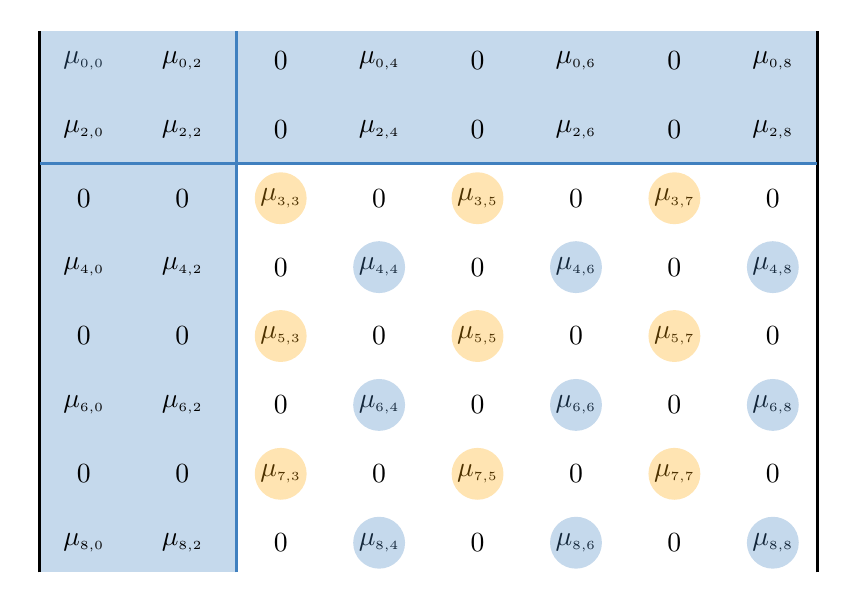}
    \caption{Moment matrix $D_9$ \eqref{eq:Dn} with the {\color{MyOra} odd} and {\color{MyBlue} even} moments are colour coded. }
    \label{fig:my_label}
\end{figure}
Note that the matrix above has a block structure with the checkerboard pattern generated by the element 
\begin{align}
    M_{i,j} = \left(\begin{array}{cc}
       {\color{MyOra}\mu_{2i-1,2j-1} } & 0 \\
        0 & {\color{MyBlue}\mu_{2i,2j}}
    \end{array} \right),\qquad i,j\geq 2.
\end{align}
The even bands are a trivial consequence of dimensions of the space of meromorphic functions in the genus 1 case.

\printbibliography
\end{document}